\documentclass[11pt,psamsfonts]{amsart}
\usepackage{amssymb,amsmath,amsthm}

\newtheorem{theorem}{Theorem}[section]
\newtheorem{corollary}[theorem]{Corollary}
\newtheorem{lemma}[theorem]{Lemma}
\newtheorem{proposition}[theorem]{Proposition}

\theoremstyle{remark}
\newtheorem{remark}[theorem]{Remark}

\theoremstyle{definition}
\newtheorem{definition}[theorem]{Definition}

\numberwithin{equation}{section}

\newcommand{\thmref}[1]{Theorem~\ref{#1}}

\newcommand{\lemref}[1]{Lemma~\ref{#1}}
\newcommand{\proref}[1]{Proposition~\ref{#1}}

\newcommand{\spec}[1]{({\sf S}$_{#1}$)}

\newcommand{\linsp}{\operatorname{span}}
\newcommand{\inv}{^{-1}}
\newcommand{\cspan}{\overline{\operatorname{span}}}

\begin{document}
\title[Partial dynamical systems]{Partial dynamical systems and
$C^*$-algebras generated by partial isometries}
\author[R. Exel]{Ruy Exel${}^1$}
\address[R. Exel]{Departamento de Matematica, Universidade de S\~{a}o Paulo,
Rua do Mat\~{a}o 1010,
05508-900 S\~{a}o Paulo, BRAZIL}
\thanks{${}^1$) Research partially supported by CNPq, Brazil.} 
\email{exel@ime.usp.br}
\author[M. Laca]{Marcelo Laca${}^2$}
\address[M. Laca]{Department of Mathematics, 
University of Newcastle, New South Wales 2308, AUSTRALIA}
\email{marcelo@math.newcastle.edu.au}
\thanks{${}^2$) Research supported by the Australian Research Council.}
\author[J. Quigg]{John Quigg${}^3$}
\address[J. Quigg]{Department of Mathematics, Arizona State University, Tempe,
Arizona 85287, U.S.A.}
\email{quigg@math.la.asu.edu}
\thanks{${}^3$) Research partially supported by
 National Science Foundation Grant No. DMS9401253}
%%%\date{November 5, 1997}
\subjclass{Primary 46L55 }
\dedicatory{{\em November 10, 1997}} 
\begin{abstract}
A collection of partial isometries whose range and initial
projections satisfy a specified set of conditions often gives rise to
a partial representation of a group. The $C^*$-algebra generated by the
partial isometries is thus a
quotient of the universal $C^*$-algebra for
partial representations of the group, from which it inherits a  crossed
product structure, of an abelian $C^*$-algebra by a partial action of the group. 
Questions of  faithfulness of representations, simplicity, and ideal structure
 of these $C^*$-algebras can then be addressed in a unified manner
{}from within the theory of partial actions. 
We do this here, focusing on two key properties of partial dynamical systems,
namely amenability and topological freeness; they are the essential ingredients
of our main results in which we characterize faithful representations,
 simplicity and the ideal structure of crossed products.
As applications we 
consider three situations involving $C^*$-algebras 
generated by partial isometries:
 partial representations of groups,
Toeplitz algebras of quasi-lattice ordered groups, and
Cuntz-Krieger algebras. These $C^*$-algebras
share a crossed product structure which we give
here explicitly and which we use to study them
 in terms of the underlying partial actions.
\end{abstract}
\maketitle
\section*{Introduction}
In this paper we develop tools to analyze partial dynamical systems
and use them in our general approach to 
$C^*$-algebras generated by partial isometries.
We realize the universal $C^*$-algebra for
partial representations of a group
subject to relations as the crossed
product by
a partial action of the group on a commutative $C^*$-algebra.
A key feature of our method is the explicit description of the
spectrum of this commutative $C^*$-algebra in 
terms of the specified relations. 
Our work builds 
upon \cite{exe-pag}, where a
certain crossed product
is shown to be universal for the partial representations of a group.

We begin by reviewing the definition and basic construction of
 crossed products by partial actions in Section \ref{CP-by-PA};
 we also establish the one-to-one correspondence between 
covariant representations of a partial dynamical system 
and representations of the associated crossed product. 

In Section \ref{topfree-sect} we 
adapt the 
notion of topological freeness for group actions \cite{arc-spi}
to the context of
partial actions on abelian $C^*$-algebras. The main technical result 
is 
\thmref{thm-fax}, where we show that for a locally compact 
Hausdorff space $X$ 
the ideals 
in the reduced crossed product of $C_0(X)$
by a topologically free partial action of a discrete group $G$   
 necessarily  intersect $C_0(X)$  nontrivially;
 hence a representation of the
reduced crossed product is faithful if and only if it is faithful 
on $C_0(X)$. This leads to a sufficient
condition for simplicity of the reduced crossed product
in Corollary \ref{simple}.

In Section \ref{amen-sect} we consider invariant ideals of 
a partial action and the ideals they generate in the crossed
product. In \proref{epimorhi-prop} we give a general short exact sequence
relating an invariant ideal of a partial action, the corresponding ideal of the
crossed product, and the crossed product by the quotient partial action. 
After discussing the approximation property 
introduced in \cite{exe-afb}, which implies amenability of a partial action 
and hence equality of the full and reduced crossed products,
we prove \thmref{ideals}, the main result of the section.
Specifically, the result is that
if a partial action has the approximation property
and is topologically free on closed invariant 
subsets, then the ideals of the crossed product are in 
one-to-one correspondence with the invariant ideals, and hence with 
invariant open  subsets under the partial action.

We begin Section \ref{spectrum-sect} by introducing a class of
partial dynamical systems arising from  partial representations 
whose range projections satisfy a given set of relations.
In \proref{spec-prop} 
we describe the spectrum associated to the relations 
and give a canonical partial action of the group on this spectrum. 
The resulting crossed product has a universal property with respect to
partial representations of the group satisfying the relations,
this is proved in \thmref{thm-spectrum}, the main result of the section. 
It is through this that the results of the first three sections
become available to study the $C^*$-algebras
generated by partial representations subject to relations.

In the final three sections we apply the main results
to three concrete situations. In Section~\ref{no-relations}
we show that the partial dynamical system canonically associated to a
discrete group in \cite{exe-pag} is topologically free if and only if the group
is infinite. Since the reduced partial 
$C^*$-algebra of such a  group
is a crossed product by a partial action,
 we are able to characterize its faithful representations.
In Section~\ref{nica-rels} we realize
the Toeplitz $C^*$-algebras associated by Nica in \cite{nica} to 
quasi-lattice ordered groups as crossed products by partial actions.
We show that the corresponding partial dynamical systems 
are topologically free, and from this we 
strengthen a result from \cite{quasilat} by showing
that a representation of the generalized Toeplitz algebra 
is faithful if and only if it is faithful on the diagonal.

In Section \ref{ck-section} we realize the Cuntz-Krieger algebra $\mathcal O_A$
associated to a $\{0,1\}$-valued $n\times n$ matrix $A$ 
as a crossed product by a partial action of the free group on $n$ generators.
We relate the properties of the spectrum of the Cuntz-Krieger relations
to properties of the matrix $A$.  In particular, we
show that the condition (I) of Cuntz and Krieger \cite{cun-kri}
corresponds to topological freeness,
and this enables us to deduce the Cuntz-Krieger uniqueness theorem from 
the results of Section~\ref{topfree-sect} and an amenability result 
{}from \cite{exe-afb}. Our description of the crossed product
realizations of the Toeplitz and the Cuntz-Krieger algebras will be
explicit, as opposed to the indirect method of \cite{qui-rae}.

{\em Acknowledgement.}  The authors would like to thank Iain Raeburn for
several helpful discussions and for his active role at the early stages of
this work.  Part of this research was conducted while the first and third
authors visited the University of Newcastle sponsored by R.M.C. visitor
grants, and they thank their hosts Marcelo Laca and Iain Raeburn for their
hospitality.

\section{Crossed products by partial actions.}\label{CP-by-PA}
Let  $\alpha$ be a partial action of the discrete group $G$ on
the  $C^*$-algebra  $A$ in the sense of \cite{ex-pa,mac,exe-pag,qui-rae}. 
We denote by $D_t$ the range of $\alpha_t$ for each $t\in G$,
and we say that the triple $(A,G,\alpha)$ is a partial dynamical system.
There are two $C^*$-algebras associated with a partial dynamical system:
the full crossed product and the reduced crossed product, cf. \cite{ex-pa,mac}.
These are defined, in analogy with the crossed products of group actions,
as certain $C^*$-completions of the convolution $*$-algebra of $A$-valued 
$\ell^1$-functions on the group. As such, they  contain the collection 
of finite sums $\{ \sum_t a_t \delta_t : a_t \in D_{t}\}$ as a dense
$^*$-subalgebra. It is also possible to view the full crossed product as a
universal $C^*$-algebra  for covariant representations 
as in \cite{qui-rae}. Since we will exploit this point of view, 
we briefly review some definitions and  basic facts. 

\begin{definition} \cite[Definition 6.2]{exe-pag}\label{def-p-rep}
A partial representation of a group $G$
on a Hilbert space $H$ is a map  $u$ from $G$
into the bounded linear operators on $H$ such that
\begin{enumerate}
\item $u(e) = 1$,
\item $u(t\inv) = u(t)^*$ and
\item $u(s) u(t) u(t\inv) = u(st) u(t\inv)$ for $s,t \in G$. 
\end{enumerate}
\end{definition}
Note that these conditions imply that the $u_t$ are 
partial isometries on $H$.
An equivalent definition which is sometimes easier to verify is given in
\cite[Definition 1.7]{qui-rae}, where one only requires that the $u_t$ 
be partial isometries with commuting range projections, satisfying
$u(e) u(e)^* = 1$, $ u(s)^* u(s) = u(s\inv) u(s\inv)^*$,
and that $ u(st) $ extends $u(s) u(t)$ in the sense of \cite[Lemma 1.6]{qui-rae}.
The equivalence is proved in \cite[Lemma 1.8]{qui-rae}.

We will use the following definition of covariant representations
of partial actions. 
\begin{definition} 
A {\em covariant representation} of 
the partial dynamical system $(A,G,\alpha)$ on a Hilbert space $H$ is a
pair $(\pi, u)$ in which $\pi$ is a nondegenerate representation
of $A$ on $H$ and $u$ is a partial representation of $G$ on $H$
 such that for each $t\in G$  we have that $u_t u_t^*$ is the
projection onto the subspace $\cspan{\pi(D_t)H}$  and 
  $$
 \pi(\alpha_t(a)) =  u_t\pi(a)u_{t\inv }, \qquad a\in D_{t\inv }.
  $$
\end{definition}

As in the case of actions of groups, covariant representations of a partial
action correspond to representations of the associated crossed product.  This
correspondence was first proved in \cite[Propositions 5.5 and 5.6]{ex-pa} in
the case of a single partial automorphism, and it was generalized to partial
actions of discrete groups in \cite[Propositions 2.7 and 2.8]{mac}.  Our
definition of covariant representations is slightly different from those of
\cite{ex-pa,mac}, but this poses no problem because the various definitions
have been shown to be equivalent \cite[Remark 1.12]{qui-rae}, see also
\cite[Section 3]{qui-rae}.  We state the results in the following Proposition
and Theorem, and include proofs based on our Definition of covariant
representations for the convenience of the reader.

\begin{proposition}
\label{proposuit}
 Let $(\pi,u)$ be a covariant representation of $(A,G,\alpha)$ on
$H$.  Then there exists a 
\textup(necessarily unique\textup)
representation, denoted
$\pi\times u$, of $A\rtimes_\alpha G$ on $H$, such that
  $$
  (\pi\times u)(a\delta_t) = \pi(a)u_t,
  $$
  for all $t$ in $G$ and all $a$ in $D_{t\inv }$.
\end{proposition}

\begin{proof}  Let $\ell^1(G,A)$ be the $\ell^1$-algebra
associated to the partial dynamical system $(A,G,\alpha)$.
If $b \in\ell^1(G,A)$, i.e., if
  $
  b = \sum_{t\in G} a_t \delta_t,
  $
  where each $a_t\in D_t$ and $\|b\|_1 = \sum_{t\in G}\|a_t\| <
\infty$, then put
  $$
  \rho(b) = \sum_{t\in G} \pi(a_t) u_t.
  $$
  Clearly, $\rho$ is a bounded linear map under the $\ell^1$-norm.  We
 claim that $\rho$ is a representation of $\ell^1(G,A)$. To prove this claim suppose
 that $a\in D_t$ and $b\in D_s$.  Then
  $$
  \rho(a\delta_t) \rho(b\delta_s) =
  \pi(a) u_t \pi(b) u_s =
  u_t u_{t\inv } \pi(a) u_t \pi(b) u_s,
  $$
  because $u_t u_{t\inv }$ is the orthogonal projection onto $H_t$,
which contains the range of $\pi(a)$.  Also, $u_{t\inv } u_t$, being
the projection onto $H_{t\inv }$, will commute with any $\pi(b)$,
since the latter leaves $H_{t\inv }$ invariant.  Therefore, the above
equals
\begin{align*}
  u_t u_{t\inv } \pi(a) u_t u_{t\inv } u_t \pi(b) u_s &=
  u_t u_{t\inv } \pi(a) u_t \pi(b) u_{t\inv } u_t u_s \\&=
  u_t \pi(\alpha_{t\inv }(a)) \pi(b) u_{t\inv } u_{ts} \\&=
  \pi( \alpha_t( \alpha_{t\inv }(a)b ) ) u_{ts}.
\end{align*}
  By definition, the product in $A\rtimes_\alpha G$ is given by
  $$
  (a\delta_t) (b\delta_s) = \alpha_t( \alpha_{t\inv }(a)b )
\delta_{ts},
  $$
hence $\rho$ is multiplicative.  We leave for the reader to
complete the proof of the claim by
verifying that $\rho$ also preserves the star operation.
 Since the crossed product is the enveloping
 $C^*$-algebra of $\ell^1(G,A)$, the *-homomorphism $\rho$ extends to
the desired representation of $A\rtimes_\alpha G$.  
\end{proof}

\begin{theorem}\label{theosuit}
Let $\alpha$ be a partial action of the group $G$ on the $C^*$-algebra $A$,
and let $H$ be a Hilbert space. Then the map
  $$
  (\pi,u) \mapsto \pi\times u
  $$
is a one-to-one correspondence between covariant representations of
$(A,G,\alpha)$ on $H$ and non-degenerate representations of
$A\rtimes_\alpha G$ on $H$.
\end{theorem}

\begin{proof}
  Let $\rho$ be a non-degenerate
representation of $A\rtimes_\alpha G$ on $H$. Identify $A$ with its
isomorphic copy $A\delta_e$ within $A\rtimes_\alpha G$ and denote by $\pi$
the restriction of $\rho$ to $A$.

For each $t\in G$ let $H_{t\inv }$ be the closure of the linear
space consisting of the vectors of the form
  $\eta = \sum_{i=1}^n \pi(a_i)\xi_i$,
  with $a_i \in D_{t\inv }$ and $\xi_i \in H$.  For
each such $\eta$ define
 \begin{equation}\label{dagger}
  u_t(\eta) = \sum_{i=1}^n \rho(\alpha_t(a_i)\delta_t)\xi_i.
  %%%\eqno{(\dagger)}
\end{equation}
  We claim that $u_t$ extends to an isometry from
  $H_{t\inv }$ onto $H_t$.  To see this let us study the norm of the
right hand side of (\ref{dagger}).  We have
  $$
  \left\Vert
  \sum_{i=1}^n \rho(\alpha_t(a_i)\delta_t)\xi_i
  \right\Vert^2 =
  \left\langle
  \sum_{i,j=1}^n
  \rho\big((\alpha_t(a_j)\delta_t)^* \alpha_t(a_i)\delta_t\big)
\xi_i,\xi_j
  \right\rangle.
  $$
 It is easy to see that
  $(\alpha_t(a_j)\delta_t)^* \alpha_t(a_i)\delta_t =
a_j^*a_i\delta_e$, and hence the above equals
  $$
  \left\langle
  \sum_{i,j=1}^n
  \pi(a_j^* a_i) \xi_i,\xi_j
  \right\rangle =
  \left\langle
  \sum_{i=1}^n \pi(a_i) \xi_i, \sum_{j=1}^n \pi(a_j) \xi_j
  \right\rangle =
  \|\eta\|^2.
  $$
It follows easily that 
 $u_t$ is a well-defined linear isometry on its domain of definition
$\linsp\{\pi(a)\xi: a \in D_{t\inv }, \xi \in H \}$. 
Hence  $u_t$  extends to an isometry,  also denoted by $u_t$, on 
$H_{t\inv } = \cspan\{\pi(a)\xi: a \in D_{t\inv }, \xi \in H \}$.
Next we show that the range of $u_t$ is precisely $H_t$.

If $\{e_\lambda\}$ is an approximate unit for $D_t$, we have,
for all $a\in D_{t\inv }$, and $\xi\in H$, that
  $$
  \rho(\alpha_t(a)\delta_t)\xi =
  \lim_\lambda \rho(e_\lambda \alpha_t(a)\delta_t)\xi =
  \lim_\lambda \pi(e_\lambda) \rho (\alpha_t(a)\delta_t)\xi,
  $$
  which shows that $\rho(\alpha_t(a)\delta_t)\xi\in H_t$ and hence
 that the range of $u_t$ is a subset of $H_t$.  Conversely, given
$a\in D_t$, observe that
  $$
  a =
  \lim_\lambda a e_\lambda =
  \lim_\lambda (a\delta_t)(\alpha_{t\inv }(e_\lambda)\delta_{t\inv }),
  $$
  and hence, for all $\xi\in H$, we have
  $$
  \pi(a)\xi =
  \lim_\lambda \rho(a\delta_t)
    \rho(\alpha_{t\inv }(e_\lambda)\delta_{t\inv })\xi =
  \lim_\lambda u_t \Big(
    \pi(\alpha_{t\inv }(a))
    \rho(\alpha_{t\inv }(e_\lambda)\delta_{t\inv })\xi
    \Big),
  $$
  proving that $H_t$ is contained in the range of $u_t$.

Consider next the extension of $u_t$ to all of $H$ defined by
setting $u_t = 0$ on the orthogonal complement of $H_{t\inv }$.  In
this way, $u_t$ becomes a partial isometry on $H$, and it is clear that its
range projection, namely $u_t u_t^*$, coincides with the orthogonal
projection onto $H_t$.

We claim that the map $t\mapsto u_t$ is a partial
representation of $G$.

We prove $u_{t\inv } = u_t^*$ first.  For
this purpose, let $a\in D_{t\inv }$ and $\xi\in H$.  We then have that
  $$
  u_{t\inv } u_t (\pi(a)\xi) =
  u_{t\inv } \rho (\alpha_t(a)\delta_t)\xi.
  $$
  Choose an approximate unit $\{e_\lambda\}$ for $D_t$, so that the
above equals
\begin{align*}
  \lim_\lambda u_{t\inv } \rho (\alpha_t(a) e_\lambda \delta_t)\xi &=
  \lim_\lambda u_{t\inv } \pi (\alpha_t(a)) \rho(e_\lambda
\delta_t)\xi \\&=
  \lim_\lambda \rho (a\delta_{t\inv }) \rho(e_\lambda \delta_t)\xi =
  \lim_\lambda \rho (a\delta_{t\inv } e_\lambda \delta_t)\xi \\&=
  \lim_\lambda \pi (a \alpha_{t\inv }(e_\lambda))\xi =
  \pi (a)\xi,
\end{align*}
  where the last equality holds because
  $\{\alpha_{t\inv }(e_\lambda)\}$ is an approximate unit for
$D_{t\inv }$.
  This shows that
  $ u_{t\inv } u_t $ coincides with the orthogonal projection onto
  $H_{t\inv }$, which is precisely the initial space of $u_t$.  In
other words
  $$
  u_{t\inv } u_t = u_t^* u_t.
  $$
  Is is evident that the initial space of $u_{t\inv }$ is $H_t$, which
is also the final space of $u_t$, that is,
  $$u_{t\inv }^* u_{t\inv } = u_t u_t^*.$$
  So, employing the last two formulas we obtain
  $$
  u_{t\inv } =
  u_{t\inv } u_{t\inv }^* u_{t\inv } =
  u_{t\inv } u_t u_t^* =
  u_t^* u_t u_t^* =
  u_t^*.
  $$

Since $\rho$ was supposed non-degenerate, it follows that $\pi$ is
also non-degenerate (because any approximate unit for $A$ is also an
approximate unit for $A\rtimes_\alpha G$).  Therefore $H_e =\pi(A)H =
H$ and it is not hard to see that $u_e$ is the identity operator on
$H$.

The last axiom of partial representations to be checked is 
  \begin{equation}\label{diamondsuit}
  u_t u_s u_{s\inv } = u_{ts} u_{s\inv },
\end{equation}
  which we now set out to prove.  A crucial observation in this
respect is that for each $t$ in $G$, the orthogonal projection onto
$H_t$, which we shall henceforth denote by $p_t$, lies in the center
of the von Neumann subalgebra of $B(H)$ generated by $\pi(A)$.  In
particular this implies that the projections $p_t = u_t u_{t\inv }$  
commute with each other.

 Since the left hand side of (\ref{diamondsuit}) satisfies
  $$
  u_t u_s u_{s\inv } =
  u_t u_{t\inv } u_t u_s u_{s\inv } =
  u_t p_{t\inv } p_s,
  $$
  we see that it vanishes on the orthogonal complement of
  $H_{t\inv } \cap H_s$. Incidentally, it is not hard to see that
  $H_{t\inv } \cap H_s = \pi(D_{t\inv }\cap D_s)H$.

Let $\xi \in \pi(D_{t\inv }\cap D_s)H$ and suppose that $\xi$ has the
form $\xi=\pi(ab)\eta$, where $a$ and $b$ are in $D_{t\inv }\cap D_s$
and $\eta \in H$.  It is clear that the closed linear span of the
set of all such $\xi$'s is equal to
  $H_{t\inv } \cap H_s$.  We have
\begin{align*}
  u_{ts}u_{s\inv }\xi &=
  u_{ts}u_{s\inv } \pi(ab)\eta =
  u_{ts} \rho(\alpha_{s\inv }(ab)\delta_{s\inv })\eta \\&=
  u_{ts} \pi(\alpha_{s\inv }(a))
    \rho(\alpha_{s\inv }(b)\delta_{s\inv })\eta \\&=
  \rho(\alpha_{ts}(\alpha_{s\inv }(a))\delta_{ts})
  \rho(\alpha_{s\inv }(b)\delta_{s\inv })\eta \\&=
  \rho(\alpha_t(a)\delta_{ts} \alpha_{s\inv }(b)\delta_{s\inv })\eta =
  \rho(\alpha_t(a)\alpha_t(b)\delta_t)\eta =
  u_t \xi.
\end{align*}
  This shows that $u_{ts}u_{s\inv }$ coincides with $u_t$ on
$H_{t\inv } \cap H_s$.  This can be expressed by saying that
  $
  u_{ts}u_{s\inv } p_{t\inv } p_s =
  u_t p_{t\inv } p_s,
  $
  which is easily seen to imply that
  \begin{equation}\label{clubsuit}
  u_{ts}u_{s\inv } p_{t\inv } =
  u_t p_s.
  %%\eqno{(\clubsuit)}
 \end{equation}
 Taking adjoints, we get
 \begin{equation}
  p_s u_{t\inv } =
  p_{t\inv } u_s u_{(ts)\inv }, 
  \end{equation}
and, replacing $s$ by  $(ts)\inv $ and $t$ by $s$, we  have 
  \begin{equation}\label{spadesuit}
  p_{(ts)\inv} u_{s\inv} = p_{s\inv}   u_{(ts)\inv} u_t
  \end{equation}
  Next, observe that the right hand side of (\ref{diamondsuit}) satisfies
  $$
  u_{ts} u_{s\inv } =
  u_{ts} p_{(ts)\inv } u_{s\inv },
  $$
  hence, using (\ref{spadesuit}) on the right-hand-side, 
  $$
  u_{ts} u_{s\inv } =
  u_{ts} p_{s\inv } u_{(ts)\inv } u_t.
  $$
  The precise form of this last identity is not absolutely crucial,
except for the fact that $u_t$ appears as the last factor in the right
hand side.  Since $u_t = u_t p_{t\inv }$, we conclude that $u_{ts}
u_{s\inv }$ is likewise not affected by right multiplication by
$p_{t\inv }$.  Using (\ref{clubsuit}), we see that
  $$
  u_{ts} u_{s\inv } =
  u_{ts} u_{s\inv } p_{t\inv } =
  u_t p_s =
  u_t u_s u_{s\inv },
  $$
  proving (\ref{diamondsuit}).

Thus $u$ is a partial representation of $G$ on $H$.  To see that it satisfies
the covariance equation
\begin{equation}\label{covarsuit}
  u_t\pi(a)u_{t\inv } = \pi(\alpha_t(a)), \qquad a\in D_{t\inv},
\end{equation}
  observe that both sides of this equation represent operators
vanishing on the orthogonal complement of $H_t$. Thus it suffices to
prove that they coincide when applied to vectors of the form
$\pi(b)\xi$, where $b\in D_t$.  We have
\begin{align*}
  u_t\pi(a) u_{t\inv } \pi(b)\xi &=
  \rho(\alpha_t(a)\delta_t)
    \rho(\alpha_{t\inv }(b)\delta_{t\inv })\xi \\&=
  \rho(\alpha_t(a)b)\xi =
  \pi(\alpha_t(a)) \pi(b)\xi,
\end{align*}
finishing the proof of (\ref{covarsuit}).

We claim next that the representation $\pi\times u$ associated to the
covariant representation $(\pi,u)$ coincides with $\rho$.  
To prove this suppose  $\eta=\pi(b)\xi$, with 
$b \in  A$, and recall that $\pi$ is non-degenerate, so the set of these
$\eta$'s is dense in $H$. Then for all $a\in D_t$
\begin{align*}
  \rho(a\delta_t) \eta &=
  \rho(a\delta_t) \pi(b) \xi =
  \rho(a\delta_t b\delta_e) \xi =
  \rho( \alpha_t(\alpha_{t\inv }(a)b)\delta_t) \xi \\&=
  \lim_\lambda \rho( \alpha_t(\alpha_{t\inv }(a) e_\lambda b)\delta_t)
\xi,
\end{align*}
where $\{e_\lambda\}$ is an approximate unit for $D_{t\inv }$.  The
above then equals
\begin{align*}
  \lim_\lambda \rho( a\alpha_t(e_\lambda b)\delta_t) \xi &=
  \lim_\lambda \pi(a) \rho(\alpha_t(e_\lambda b)\delta_t) \xi =
  \lim_\lambda \pi(a) u_t \pi(e_\lambda b) \xi \\&=
  \lim_\lambda \pi(a) u_t \pi(e_\lambda) \pi(b) \xi =
  \pi(a) u_t \pi(b) \xi,
\end{align*}
  where the last step follows from the fact that $\pi(e_\lambda)$
converges, in the weak operator topology, to the projection onto
$H_{t\inv }$, which is the initial space of $u_t$.
  This shows that $\pi\times u = \rho$.

To conclude we must prove that $(\pi,u)$ is the unique covariant
representation such that $\pi\times u = \rho$.  Since we want
$\rho(a\delta_e) = \pi(a) u_e$, then clearly $\pi$ is uniquely
determined.  Assuming that $(\pi,v)$ is another 
covariant representation such that $\pi \times v = \rho$, 
we have, for $a,b\in
D_{t\inv }$ and $\xi\in H$, that
\begin{align*}
  v_t \pi(ab)\xi &=
  v_t \pi(a) v_{t\inv } v_t \pi(b)\xi =
  \pi(\alpha_t(a)) v_t \pi(b)\xi \\&=
  \rho (\alpha_t(a)\delta_t) \pi(b)\xi =
  u_t \pi(a) \pi(b)\xi =
  u_t \pi(ab)\xi.
\end{align*}
  Since the collection of vectors of the form $\pi(ab)\xi$, as above,
is dense in $H_{t\inv }$, we conclude that $v_t=u_t$.
\end{proof}
 
\section{Topologically free partial actions.}\label{topfree-sect}
We will mostly be concerned with partial actions 
arising from partial homeomorphisms of a locally compact space
$X$, so that for every $t\in G$ there is an open subset $U_t$ of $X$ 
and a homeomorphism $\theta_t : U_{t\inv} \to U_{t}$ such that $\theta_{st}$
extends  $\theta_s \theta_t$ \cite{exe-pag,mac,ex-tpa}.
 The  partial action $\alpha$ of $G$ on
$C_0(X)$ corresponding to $\theta$ is given by
 $$\alpha_t(f)(x) := f(\theta_{t\inv}(x)), \qquad f \in C_0(U_{t\inv}).$$
So, here the ideals are $D_t=C_0(U_t)$.
 We will talk about the partial action at either 
the topological or the $C^*$-algebraic level, according to convenience.

\begin{definition}[cf. \cite{arc-spi}]
The partial action $\theta$ is {\em topologically free}
if for every $t \in G \setminus \{e\}$ the set
$F_t := \{x\in U_{t\inv} : \theta_{t}(x) = x \}$
has empty interior.
\end{definition}
We point out that although the set $F_t$ need not be closed in $X$,
it is relatively closed in the domain $U_{t\inv}$ of $\theta_t$.
A standard argument gives the following equivalent 
version of topological freeness which is more
appropriate for our purposes. 

\begin{lemma}\label{nowh-dense}
The partial action $\theta$ on $X$ is topologically free if and only if 
 for every finite subset  $\{t_1, t_2, \ldots, t_n\}$ of $G \setminus \{e\}$, 
the set $\bigcup_{i = 1 }^n F_{t_i}$ has empty interior.
\end{lemma} 
\begin{proof}
It suffices to show that for every $t \in  G \setminus \{e\}$, the fixed point set
$F_t$ is nowhere dense (i.e., its closure in $X$ has empty interior), and then use 
the fact that a finite union of nowhere dense sets is nowhere dense.

Since  $F_t$ is closed relative to $U_t$ we can write $F_t = C \cap U_t$ with $C$
closed in $X$. Suppose $V \subset \overline{F_t}$ is open. 
Since the set $V\cap U_t$ is contained in $C\cap U_t = F_t$, it must be empty, by
the assumption of topological freeness. Thus $V$ and $U_t$ are disjoint open sets, 
so each one is disjoint from the other's closure. But $V \subset \overline{C \cap
U_t} \subset C \cap \overline{U_t} \subset \overline{U_t}$, 
so $V$ itself is empty and hence $F_t$ is nowhere dense. 
\end{proof}
\begin{lemma}\label{hlemma}
Let $t \in G\setminus \{e\}$, $f\in D_t$,
and $x_0 \notin F_t$.
For every $\epsilon > 0 $ there exists $h \in C_0(X)$ such that
\begin{enumerate}
\item  $h(x_0) = 1$,
\smallskip
\item $\| h  ( f \delta_t) h \| \leq \epsilon $, and
\smallskip
\item $ 0 \leq h \leq 1$.
\smallskip
\end{enumerate}
\end{lemma}
\begin{proof}
We separate the proof into two cases according to $x_0$
being in the domain $U_t$ of $\theta_{t\inv}$ or not.
If $x_0 \notin U_t$, let $K:= \{x\in U_t: |f(x)| \geq \epsilon\}$.
Then $K$ is a compact subset of $D_t$ and $x_0 \notin K$, so there is
$h\in C_0(X)$ such that $0 \leq h \leq 1$, $h(x_0) = 1$ and $h(K) = 0$.
Since $f$ is bounded by $\epsilon$ off $K$, it follows that
$\| h f \| \leq \epsilon$, so  (ii) holds too.

If $x_0 \in U_t$ then $\theta_{t\inv}(x_0)$ is defined and not equal to $x_0$.
Take disjoint open sets $V_1$ and $V_2$ such that 
$x_0 \in V_1$ and $\theta_{t\inv}(x_0) \in V_2$.
We may assume that $V_1 \subset U_t$ and $V_2 \subset U_{t\inv}$.

Letting $V: = V_1 \cap \theta_t(V_2)$, we have that $x_0 \in V \subset V_1$ and 
$\theta_{t\inv} (V) \subset V_2$, from which it follows that
$\theta_{t\inv} (V) \cap V = \emptyset$. Take now $h \in C_0(X)$ such that
$0 \leq h \leq 1$, $h(x_0) = 1$ and $h(X\setminus V) = 0$.
It remains to prove that $h$ satisfies (ii). In fact, the product
$h f \delta_t h = \alpha_t( \alpha_{t\inv} (h f) h ) \delta_t$ 
vanishes because the support of $\alpha_{t\inv} (h f)$ is contained in 
$\theta_{t\inv}(V)$ and the support of $h$ is in $V$.
\end{proof}

The reduced crossed product  associated to a partial dynamical system
in \cite[\S 3]{mac} can also be obtained as the reduced cross-sectional 
algebra of the Fell bundle determined by  the partial 
action \cite[Definition 2.3]{exe-afb}.

This reduced crossed product is a topologically graded algebra and 
the conditional expectation, denoted by $E_r$, from $C_0(X) \rtimes_r G$
onto $C_0(X)$ is a faithful positive map \cite[Proposition 2.12]{exe-afb}
(see also \cite[Corollary 3.9 and Lemma 1.4]{qui-dis} and \cite[Corollary
3.8]{qui-rae}).

\begin{proposition}\label{hprop}
If  $(C_0(X), G, \alpha)$ is a topologically free partial action then for every 
$c \in C_0(X) \rtimes_r G$ and every $\epsilon > 0$ there exists $h \in C_0(X)$ 
such that:
\begin{enumerate}
\item   $\| h E_r(c) h \| \geq \| E_r(c) \| - \epsilon$,
\smallskip
\item  $\| h E_r(c) h  - h c h \| \leq \epsilon $, and
\smallskip
\item  $ 0 \leq h \leq 1$.
\smallskip
\end{enumerate}
\end{proposition}

\begin{proof}
Assume first $c$ is a finite linear combination 
of the form $\sum_{t\in T} a_t \delta_t$, where $T$ denotes a finite subset of $G$,
in which case $E_r(c) = a_e$ (where we put $a_e = 0$ if $e \notin T$). 
Let $V = \{x \in X: |a_e(x)| > \|a_e\| - \epsilon\}$, which is
clearly open and nonempty. By \lemref{nowh-dense}
there exists $x_0 \in V$ such that $x_0  \notin  F_t$ for every 
$t \in T \setminus \{e\}$, and by \lemref{hlemma}
there exist functions $h_t$ satisfying 
$$
h_t(x_0) = 1,
\quad
\|h_t (a_t \delta_t) h_t \| \leq \frac{\epsilon}{|T|},
\quad\text{and}\quad
0 \leq h_t \leq 1.
$$
Let $h := \prod_{t\in T\setminus \{e\}} h_t$. Then (iii) is 
immediate, and (i) also holds because
 $x_0 \in  V$ so $\|h a_e h \| \geq |a_e(x_0)| > \|a_e \| - \epsilon$.
For (ii), we have
\begin{eqnarray*}
\| h a_e h - h a h \| & = &
  \| \sum_{t \in T \setminus \{e\}} h a_t \delta_t h\ \| \\
& \leq & \sum_{t \in T \setminus \{e\}} \| h a_t \delta_t h \| \\
& \leq & \sum_{t \in T \setminus \{e\}} \| h_t a_t \delta_t h_t \| \\
&  < & \epsilon.
\end{eqnarray*}

Since the elements of the form $\sum_{t \in T} a_t \delta_t$ are dense
in the crossed product and the conditional expectation $E_r$ is contractive,
a standard approximation argument gives the general case.
\end{proof}

\begin{remark}
It is easy to see that the Proposition also holds with the full crossed product
replacing the reduced one.
\end{remark}

\begin{theorem}\label{thm-fax}
Suppose $(C_0(X), G, \alpha)$ is a topologically free partial action. If $I$ is
an ideal in $C_0(X) \rtimes_r G$ with $I \cap C_0(X)  = \{0\}$, then $I = \{0\}$.
A representation of the reduced crossed product $C_0(X) \rtimes_r G$
is faithful if and only if it is faithful on $C_0(X)$.
\end{theorem}

\begin{proof}
Denote by $\pi: C_0(X) \rtimes_r G \to (C_0(X) \rtimes_r G)/ { I}$ the quotient
map, and let $a \in  { I}$ with $a \geq 0$, so that $\pi(a) = 0$.
Given $\epsilon > 0$ take $h \in C_0(X)$ satisfying conditions
(i), (ii) and (iii) of \proref{hprop}. Then 
$$
\| \pi ( h E_r(a) h) \| = \| \pi ( h (E_r(a) - a) h) \| \leq \epsilon,
$$
because $\pi(a) = 0 $. Since  $\pi$ is isometric on $C_0(X)$, because 
${ I} \cap C_0(X) = \{0\}$, it follows that 
$ \| h E_r(a) h \| \leq \epsilon$.
By \proref{hprop} (i),  $\| E_r (a)\| - \epsilon \leq \| h E_r(a) h\| $, so
$\| E_r(a) \| \leq 2 \epsilon$, and $E_r(a)$ has to vanish.
Since the conditional expectation $E_r$ is faithful on the reduced crossed
product this implies that $a = 0$ and hence that ${ I} = \{0\}$.
This proves the first assertion, which, applied to the kernel of a representation,
gives the second one.
\end{proof} 

\begin{definition}
A subset $V$ of $X$ is {\em invariant} under the partial action $\theta$
on $X$ if $\theta_s(V \cap U_{s\inv}) \subset V$ for every $s \in G$.

An ideal $J $ in $C_0(X)$ is {\em invariant} under 
the corresponding partial action $\alpha$ on $C_0(X)$ if
$\alpha_t( J \cap D_{t\inv}) \subset J$ for every $t\in G$.
\end{definition}

 It is easy to see that if $U$ is an invariant open set
then the associated ideal $C_0(U)$ is invariant and, conversely,
 every invariant ideal corresponds to an invariant open set.

\begin{definition}
The partial action $\theta$ on $X$ is {\em minimal} if there are no $\theta$-invariant
open subsets of $X$ other than $\emptyset$ and $X$ or, equivalently, if the
partial action $\alpha$ on $C_0(X)$ has no nontrivial proper invariant ideals.
\end{definition}

The complement of an invariant set is invariant too, so the partial action is
minimal if and only if it has no nontrivial proper closed invariant subsets. 

\begin{corollary}\label{simple}
If a partial action is topologically free and minimal then
the associated reduced crossed product is simple.
\end{corollary}
\begin{proof}
Suppose ${J}$ is the kernel of a representation $\rho$ 
of $C_0(X) \rtimes_r G$, and write $\rho = \pi \times v$ by \thmref{theosuit}.
Then  $J \cap C_0(X) $ is an ideal in $C_0(X)$ which is invariant under $\alpha$ 
because for every $f \in  J \cap D_{t\inv}$, we have, by covariance, that 
$\pi ( \alpha_t( f) ) = v_t \pi (f) v_t^* = 0$, and hence that $\alpha_t(f) \in J$.

By assumption $\alpha $ is minimal, so either
$J\cap C_0(X) = C_0(X)$, in which case $\pi = 0$, hence $\rho = 0$
by \proref{proposuit}, or else   $J\cap C_0(X)=\{0\}$, in which case
the representation $\pi$ is faithful by \thmref{thm-fax}.
This proves that the crossed product is simple.
\end{proof} 

\section{Invariant ideals and the approximation property.}\label{amen-sect}

Let $\alpha$ be a partial action on the $C^*$-algebra $A$.
For each invariant ideal $I$ of $A$ there is a restriction of $\alpha$ to a
partial action on $I$, with ideals $D_t\cap I$ as domains of the restricted
partial automorphisms, and there is also a quotient partial action
$\dot{\alpha}_t$ of $G$ on $A/I$, defined  by composition with the quotient
map $a\in A \mapsto  a + I \in A/I$:  the domain of $\dot{\alpha}_t$
is the ideal $\dot{D}_{t\inv} : = \{ a+I \in A/I: a \in D_{t\inv}\}$
and $\dot{\alpha}_t(a+I) =   \alpha_t(a) + I $. 

We will show that the quotient of the crossed product 
$A\rtimes G$ modulo the ideal generated by $I$ is isomorphic
to the crossed product of the quotient partial action modulo $I$.
This result generalizes \cite[Proposition 3.4]{qui-dua}, which proves the case 
$G = \mathbb Z$, and extends part of \cite[Proposition 5.1]{mac}, which only
concerns ideals. The original argument, for group actions, is from
\cite[Lemma 1]{gre-smo}. We will denote by $\langle S \rangle$ the ideal
generated by a subset $S$ of a $C^*$-algebra $B$. 

\begin{proposition}\label{epimorhi-prop}
Suppose $\alpha$ is a partial action on $A$ and assume $I$ is an 
$\alpha$-invariant ideal of $A$. Then the map 
$a\delta_t \in I\rtimes G  \mapsto a \delta_t \in A\rtimes G$  extends to an
injection of $I \rtimes G $ onto the ideal $ \langle I \rangle$ generated by 
$I$ in $A \rtimes G$, and $\langle I \rangle \cap A = I$.

The map $a \delta_t \in A\rtimes G \mapsto (a+I) \delta_t \in (A/I) \rtimes G$ 
extends to a surjective homomorphism, giving the short exact sequence
$$
0 \to  I \rtimes G  \to  A \rtimes G  \to  (A/I) \rtimes G \to 0.
$$
\end{proposition} 

\begin{proof}
The assertion that $I\rtimes G$ injects as an ideal in $A\rtimes G$ is proved in
\cite[Proposition 5.1 and Corollary 5.2]{mac} and, as done there,
we identify $I\rtimes G$ with $ \cspan \{a \delta_t: a\in D_t \cap I, \, t\in G\}$;
we also identify $I$ with its canonical image $I \delta_e$ in $A\rtimes G$.
It is clear that
$\langle I\rangle $ is contained in $I \rtimes G$. To prove the reverse inclusion
it suffices to show that $a\delta_t \in \langle I\rangle$ 
for every $a \in D_t \cap I$ and $t\in G$. Assume $a \in D_t \cap I$ and let
$b_\lambda$  be an approximate unit for the ideal $D_t$. Then 
$a b_\lambda \delta_t = (a\delta_e)(b_\lambda\delta_t)  \in \langle I\rangle $
so  $a\delta_t = \lim_\lambda a b_\lambda \delta_t \in \langle I\rangle $.
This proves that $I\rtimes G = \langle I \rangle$, from which
it is obvious that $I = \langle I \rangle \cap A$.

The map $ a \delta_t \mapsto (a+I) \delta_t$ induces a ${}^*$-homomorphism from 
$\ell^1(G,A)  $ onto  $\ell^1(G,A/I)$. Since $A\rtimes G$ is the enveloping
$C^*$-algebra  of $\ell^1(G,A) $, there is  $C^*$-homomorphism  
$\phi$ of $A\rtimes G$ onto $(A/I) \rtimes G$
which sends $ a \delta_t$ to $(a+I) \delta_t$ 
for every $a \in D_t$ and every $t\in G$.
To finish the proof we need to show that $\ker \phi = \langle I \rangle$.

It is clear that $\ker \phi$ contains the ideal $ \langle I \rangle$
generated by $I$ in $A \rtimes G$.
It remains to prove that  $\ker \phi \subset \langle I \rangle$. 
Let $\pi \times u$ be a representation of  $A \rtimes G$
with kernel $\langle I \rangle$. Since the kernel of $\pi$
contains $I$, $\pi$ factors through the quotient map $A \to A/I $;
 denote by $\tilde{\pi}$ the corresponding representation of $A/I$.
The pair $(\tilde{\pi}, u)$ is covariant and determines a representation
$\tilde{\pi} \times u$ of $(A/I) \rtimes G$.
Then $\pi \times u = (\tilde{\pi} \times u ) \circ \phi $, so 
$\ker \phi \subset \ker(\pi \times u)$.
\end{proof}
\begin{remark}
We point out that, at the level of {\em reduced}
crossed products, it is always true that $I\rtimes_r G$ injects
as an ideal in $A\rtimes_r G$ \cite[Proposition 5.1]{mac},
 but whether the quotient is the reduced crossed
product $(A/I) \rtimes_r G$ is a subtler question. We refer to the
discussion at the end of \cite[\S 4]{exe-afb} for related considerations.
\end{remark}

When $A = C_0(X)$ the $\alpha$-invariant 
ideals are in one to one correspondence with
invariant open sets; the corresponding quotients are naturally identified
with the continuous functions on the complements of these invariant open sets.
Specifically, if $\alpha$ is a partial action on $C_0(X)$ and
 $U$ is an invariant open subset of $X$ 
 then $C_0(U)$ is an invariant ideal in $C_0(X)$,
and every invariant ideal is of this form.
Moreover, $C_0(X) / C_0(U) \cong C_0(\Omega)$  with $\Omega = X \setminus U$,
the quotient map being simply restriction.
The quotient partial action 
$\dot{\alpha}$ of $G$ on  $C_0(\Omega)$ is given by 
$\dot{\alpha}_t(f|_\Omega) =  \alpha_t(f)|_\Omega$  for $f\in D_{t\inv}$
(the domain of $\dot{\alpha}_t$ 
consists of the restrictions to $\Omega$ of functions in $D_{t\inv}$).

In general there may be more ideals
in a crossed product than those of the form  
$\langle I \rangle$ with $I$ an invariant ideal in $A$.
Easy examples abound even for full actions; for instance 
write $C^*(G) = \mathbb C \rtimes G$  
(with the trivial action).
If $G$ has more than one element, then 
the kernel of the trivial homomorphism $s\mapsto 1$ from $C^*(G)$ to
$\mathbb C$
is a proper nontrivial ideal which is not generated 
by an ideal in $\mathbb C$.

The quotient system $(C_0(\Omega), G, \alpha)$, obtained by 
restricting the partial action $\alpha$
to a closed invariant subset $\Omega$ of $X$,
need not be topologically free even if $(C_0(X), G, \alpha)$ is.
An easy example of this phenomenon is obtained by restricting 
the action of $\mathbb Z$ by translation on 
$C(\mathbb Z \cup \{\pm \infty\})$ to the subset $\{\pm \infty\}$.
However, we will see that if topological freeness holds on quotients
of a partial action having the approximation
property introduced in \cite{exe-afb}, then
all the ideals of the crossed product are obtained from their 
intersections with $C_0(X)$, via the map $I \mapsto \langle I \rangle$.
Before we prove this we briefly review some basic facts about
amenability and the approximation property.

A partial dynamical system $(C_0(X), G, \alpha)$
is {\em amenable} if the canonical homomorphism from the full crossed 
product to the reduced one is faithful. Amenability is
equivalent to faithfulness (as a positive map) of the conditional
expectation from the full crossed product $C_0(X) \rtimes G$
onto $C_0(X)$ 
 \cite[Proposition 4.2]{exe-afb}; it is is also equivalent 
to {\em normality} of the dual coaction in the sense of \cite{qui-fulred}
and to amenability of the associated semi-direct product Fell bundle 
\cite[Definition 2.8]{ex-tpa}. 

\begin{definition} 
The partial dynamical system $(C_0(X), G, \alpha)$ has the { \em approximation
property} if the semi-direct product bundle has the approximation property of
\cite[Definition 4.4]{exe-afb}, that is, if there exists a net $(a_i)$ of
finitely supported functions $a_i : G \to C_0(X)$ such that
$$
\sup_{i} \| \sum_{t \in G} a_i(t)^* a_i(t) \| < \infty
$$
and 
$$
\lim_{i} \sum_{t \in G} a_i(st)^* f \delta_s  a_i( t) 
= f\delta_s \qquad s\in G, f\in D_{s}.
$$
\end{definition}
This approximation property implies amenability 
\cite[Theorem 4.6]{exe-afb}, over which it has the advantage of being
inherited by graded quotients of Fell bundles
in the sense specified in the next proposition.
We do not know at present whether the approximation property 
is actually equivalent to  amenability. 

Although we will only need the following result in the
special situation of crossed products by partial actions, 
it is more convenient to formulate it
for the topologically graded algebras studied in \cite{exe-afb}.

\begin{proposition}\label{app-prop-quo}
Suppose the  $C^*$-algebra  $B$ is topologically graded over $G$, 
and assume the associated Fell
bundle has the approximation property. Let $J \subset B$ be such that
$J = \langle J\cap B_e \rangle$.
Then the quotient $B/J$ is topologically graded over $G$ and its
associated Fell bundle also has the approximation property.
\end{proposition}
\begin{proof}
Let $\pi : B \to B/J$ be the quotient map.
That $B/J$ is topologically graded is proved in \cite[Proposition 3.1]{exe-afb}.
To prove the second claim, suppose that
the $a_i$ are the approximating functions for $B$.
Then the collection of functions $t \mapsto \pi(a_i(t))$
can be used to show that the approximation property holds for $B/J$. 
\end{proof}

\begin{theorem}\label{ideals}
Let $(C_0(X), G, \alpha)$ be a partial dynamical system
which is topologically free 
on every closed invariant subset of $X$ and which 
satisfies the approximation property. Then the map
$$ U \mapsto \langle C_0(U) \rangle $$
is a lattice isomorphism of the invariant open subsets
of $X$ onto the ideals of $C_0(X) \rtimes G$.
\end{theorem}

\begin{proof} 
It is clear that the map $U \mapsto \langle C_0(U)\rangle$ maps invariant open
subsets of $X$ to ideals in $C_0(X) \rtimes G$, and that
if $U_1 \subset U_2$ then 
$\langle C_0(U_1)\rangle \subset \langle C_0(U_2)\rangle$.
Since $C_0(U) = \langle C_0(U)\rangle \cap C_0(X)$
by \proref{epimorhi-prop},
 the map is one-to-one.
Next we show that every ideal in $C_0(X) \rtimes G$ is of this 
form; this will prove that $U \mapsto \langle C_0(U)\rangle$
is an order preserving bijection, hence a lattice isomorphism.

Suppose $J$ is an ideal of $C_0(X) \rtimes G$, and let
$I := J \cap C_0(X)$. Then $I = C_0(U)$ for an open invariant subset
$U \subset X$, and it is clear that $\langle C_0(U) \rangle \subset J$;
we will show that in fact $\langle C_0(U) \rangle = J$.

The set $\Omega := X \setminus U$ is closed and invariant, 
so $\langle C_0(U) \rangle$ is 
the kernel of the homomorphism 
$$
\phi: C_0(X) \rtimes G \to  C_0(\Omega) \rtimes G
$$ 
by the preceding proposition.
Let $b \in \phi(J) \cap C_0(\Omega)$, so that $b = \phi(a)$
for some $a \in J$ and $b = \phi(a_1) $ for some $a_1 \in C_0(X)$.
Thus  $a - a_1 \in \ker \phi$,  and since 
 $\ker \phi = \langle C_0(U) \rangle \subset J$ it follows that
 $a_1$ itself is in $J$. But then $a_1 \in J \cap C_0(X)
= C_0(U)$, so $b = \phi(a_1) = 0$.
This shows that the ideal $\phi(J)$  of $C_0(\Omega) \rtimes G$
has trivial intersection with $C_0(\Omega)$.  
  
By \proref{app-prop-quo} the partial action on the quotient
$C_0(\Omega)$ satisfies the approximation property, so it
is amenable and the reduced and full crossed products coincide, 
by \cite[Proposition 4.2]{exe-afb}.

Since by assumption $\alpha$ is topologically free on $\Omega$, 
 $\phi(J)$ is trivial by  Proposition \ref{thm-fax}, and thus
$J \subset \ker \phi = \langle C_0(U) \rangle$ as required.
\end{proof} 

\section{Partial representations subject to conditions.}\label{spectrum-sect} 

We begin by reviewing some of the main ideas from \cite{exe-pag}.
Consider the compact Hausdorff space $\{0,1\}^G$, and let $e$ 
denote the identity element in $G$. The subset
$$
X_G := \{\omega \in \{0,1\}^G : e\in \omega\}
$$
is a compact Hausdorff space with the relative topology inherited from $\{0,1\}^G$. 

The sets $X_{t} := \{ \omega \in X_G: t \in \omega\}$ 
are clopen, and we define a partial homeomorphism
$\theta_t$ on $X_{t\inv}$ by  $\theta_t(\omega) = t\omega$,
 where $t\omega = \{tx: x\in \omega\}$. This gives a  partial action
$(\{X_t\}, \{ \theta_t\} )_{t\in G}$ canonically associated to the group $G$.

At the algebra level, denote by $1_t$ the characteristic function of $X_t$;
then $C(X_G)$ is the closed linear span of the projections 
$\{1_s: s\in G\}$.
The domain of the partial automorphism $\alpha_t$ is 
$C_0(X_{t\inv}) = \cspan\{ 1_s 1_{t\inv} : s \in G\}$, and $\alpha_t$ 
is determined by  $$\alpha_t(1_s 1_{t\inv} )   = 1_{ts} 1_t .$$

The crossed product 
$C(X_G) \rtimes_\alpha G$ has the following universal property:
\begin{itemize}
\item[(U1)] For every partial representation $u$ of $G$ 
there is a unique representation $\rho_u$
of $C(X_G)$ satisfying $\rho_u(1_t) = u_t u_t^*$, 
and  $(\rho_u, u)$ is a covariant representation of 
$(C(X_G), G, \alpha)$;
\item[(U2)] every 
representation of $C(X_G) \rtimes_\alpha G$ is of the form $\rho_u \times u$
with $(\rho_u, u)$ as above.
\end{itemize}
Since $C^*(\{u_t: t\in G\})$ coincides with the $C^*$-algebra generated 
by the range of $\rho \times u$,
this justifies referring to the crossed product $C(X_G) \rtimes G$ as 
the universal $C^*$-algebra for partial representations
of $G$ or, simply, as the {\em partial group algebra} of $G$ 
\cite[Definition 6.4]{exe-pag}.
 
Notice that since $X_t$ is clopen the partial isometries themselves belong to
 the crossed product and, indeed, they generate it. We will denote by
$[t]$ the partial isometry corresponding to the group element $t$
in the universal partial representation of $G$, and, by abuse of notation,
we will also write $1_t = [t][t]^*$. When $u$ is a partial representation 
of $G$ we will denote the range projections $u_t u_{t\inv} = u_t u_t^* $ by
$e^{u}(t)$  or simply by $e(t)$. Notice that the initial projection of $u_t$ is 
the range projection of $u_{t\inv}$, so we only need to mention range projections.

We are interested here in  partial representations whose range projections
satisfy a set of relations of the form
$$
 \sum_i \prod_j \lambda_{ij} e(t_{ij}) = 0,
$$
where the $\lambda_{ij}$ are scalars and the sums and products are over finite sets. 
More generally, given a collection of functions $\mathcal R $
in $C(X_G)$ we will say that
the partial representation $u$ of $G$ 
{\em satisfies the relations} $\mathcal R$ if
the representation $\rho_u$ of $C(X_G)$ obtained by extending the map
$1_t \mapsto u_t u_t^*$ vanishes on every $f \in \mathcal R$.
When the relations in $\mathcal R$ are of the form specified above, 
this amounts to saying that the generating partial isometries satisfy 
$\sum_i \prod_j \lambda_{ij}u_{t_{ij}} u_{t_{ij}}^* = 0$.

\begin{proposition}\label{spec-prop}
Let $\mathcal R$ be a collection of functions in $C(X_G)$. Then the smallest
$\alpha$-invariant \textup(closed, two-sided\textup) ideal of $C(X_G)$
containing $\mathcal R$ is the ideal, denoted $I$, generated by
the set  $\{ \alpha_t(f 1_{t\inv}) : t \in G, f \in \mathcal R \}$.
Moreover, 
\begin{equation}\label{spec-def}
\Omega_{\mathcal R} := \{ \omega \in X_G : 
f(t\inv \omega) = 0 \text{ for all } t \in \omega, f \in
\mathcal R\}
\end{equation}
is a compact invariant subset of $X_G$ such that $I=C_0(X_G\setminus
\Omega_R)$, and the quotient $C(X_G) / I$ is 
canonically isomorphic to $C(\Omega_{\mathcal R})$.
\end{proposition}

\begin{proof}
Notice first that for every $f \in C(X_G)$ the function 
$f 1_{t\inv} $ is in $D_{t\inv}$ so that it makes sense to talk about
 $\alpha_t(f 1_{t\inv})$. Moreover, 
identifying $ C(X_G)$ with its image in the crossed product 
and using covariance, we have  $\alpha_t(f 1_{t\inv})  = [t] f [t\inv]$.
Let  $I$ be the ideal generated by
 $ \{ \alpha_t(f 1_{t\inv}) : t \in G, f \in \mathcal R\}$.
Since any invariant ideal which contains $f$ must contain 
$\alpha_t(f 1_{t\inv}) = [t] f [t\inv] $, 
the smallest invariant ideal containing
$\mathcal R$ must contain $I$.  
The reverse inclusion will follow once
we show that $I$ is invariant, i.e., that
$\alpha_t(I\cap D_{t\inv}) \subset I$ for every $t\in G$.
Since  $I\cap D_{t\inv} =  1_{t\inv} I$,  we need to show that
$\alpha_t(f 1_{t\inv}) \in I$ for every $f\in I$ and $t\in G$.
For $g \in C(X_G)$, $ f\in \mathcal R$ and $s,t\in G$, we have
\begin{eqnarray*}
[s] \alpha_t(f 1_{t\inv}) g [s\inv] & =&[s] ([t] f [t\inv] g ) [s\inv]\\
& = &  [st][t\inv][t] f [t\inv] g [s\inv][s][s\inv] \\
&  =& 
 [st] f [t\inv][s\inv][s] g [s\inv] \\
& = & [st] f [(st)\inv] g'
\end{eqnarray*}
with $g' = [s] g [s\inv]$.
 Since the linear span of the elements $\alpha_t(f 1_{t\inv}) g$ is dense
in $I$ and $\alpha_t(f 1_{t\inv}) = [t] f [t\inv]$, invariance follows.

Let $U \subset X_G$
be the invariant open set such that $I = C_0(U)$. 
Then the quotient $C_0(X_G) / I$ 
is isomorphic to $C(X_G \setminus U)$ and it only remains to
prove that $U = X_G \setminus \Omega_{\mathcal R}$.

Since $\alpha_t(f 1_{t\inv}) (\omega)$ 
is equal to $f(t\inv \omega)$ when $t\in \omega$,
and $0$ otherwise, the characterization of $I$ given in the first part
implies that  $f(t\inv \omega) = 0$ for every $t\in \omega$ and $f\in \mathcal R$
 if and only if $F(\omega) = 0$ for every $F \in I$. 
This  proves that $\Omega_{\mathcal R} = X_G \setminus U$, finishing the proof.  
\end{proof}
    
\begin{definition}
The set $\Omega_{\mathcal R}$ is called 
the {\em spectrum of the relations} ${\mathcal R}$.
\end{definition}

The spectrum of a set of relations is invariant under
the partial action $\alpha $ on $C(X_G)$ so there is a quotient
 partial action (also denoted $\alpha$) on
$C(\Omega_{\mathcal R})$ obtained by restricting the partial homeomorphisms
to $\Omega_{\mathcal R}$.
The restricted partial homeomorphisms have compact 
open (relative to $\Omega_{\mathcal R}$)
sets as domains and ranges so for each
group element $t$ the partial isometry $v_t = (v_t v_t^*) v_t$
belongs to the crossed product. We will show that the  crossed product
$C(\Omega_{\mathcal R}) \rtimes G$ has a universal property
with respect to partial representations of $G$ subject to the relations $\mathcal R$.
 
\begin{definition}
Suppose $G$ is a group and let $\mathcal R \subset C_0(X_G)$ 
be a set of relations. A partial representation $v$ of $G$ is 
{\em universal for the relations $\mathcal R$} if  
\begin{enumerate}
\item $v$ satisfies $\mathcal R$, i.e., $\rho_v(\mathcal R) = \{0\}$, and
\item for
every partial representation $V$ of $G$ satisfying $\mathcal R$
the map $v_t \mapsto V_t$ extends to a $C^*$-algebra homomorphism
{}from $C^*(\{ v_t: t \in G\})$ onto $C^*(\{ V_t: t \in G\})$.
\end{enumerate}

The $C^*$-algebra generated by a universal partial representation for 
$\mathcal R$ (which is clearly unique up to canonical isomorphism) will 
be called the {\em universal $C^*$-algebra for partial representations of 
$G$ subject to the relations $\mathcal R$} and denoted $C^*_p(G;\mathcal R)$. 
\end{definition}

The existence of universal representations subject to relations,
and of the universal $C^*$-algebra for $\mathcal R$,
could be derived from an abstract argument, cf. \cite{bla}; 
we choose instead to give a concrete realization as a crossed product.

\begin{theorem}\label{thm-spectrum}
Suppose $\mathcal R$ is a collection of relations
in $C(X_G)$  with spectrum 
$$
\Omega_{\mathcal R}  = \{ \omega \in X_G : f (t\inv\omega) 
= 0 \text{ for all } t \in \omega ,f \in \mathcal R\}.
$$

\begin{enumerate}
\item  If $\rho\times V$ is a representation of 
$C(\Omega_{\mathcal R}) \rtimes G$ then $V$
is a partial representation of $G$ 
satisfying the relations $\mathcal R$.

\item Conversely, if $V$ is
a partial representation of $G$ satisfying the relations $\mathcal R$,
then  $1_{t} \mapsto V(t) V(t)^*$ extends uniquely to a representation $\rho_V$ of
$C(\Omega_{\mathcal R})$, and the pair $(\rho_V, V)$ is covariant.

\item  The universal $C^*$-algebra $C^*_p(G;\mathcal R)$ for
partial representations of $G$ subject to the relations ${\mathcal R}$ exists and
is canonically isomorphic to  $C(\Omega_{\mathcal R}) \rtimes G$.
\end{enumerate}
\end{theorem}

\begin{proof}
(i)  holds because the
range projections of the partial isometries $v_t$
 are in $C(\Omega_{\mathcal R})$, which was defined
precisely so that the 
relations $\mathcal R$
be satisfied.

Next we prove (ii). If $V$ is a partial representation satisfying the relations, then
the representation of $C(X_G)$ determined by the range projections $1_t \mapsto V_t V_t^*$ 
factors through $C(\Omega_{\mathcal R})$, and hence  gives a covariant representation of
 $(C(\Omega_{\mathcal R}), G, \alpha)$. 

 By (i) and (ii) there is a bijection between 
 partial representations satisfying the relations and  covariant
representations of $(C(\Omega_{\mathcal R}), G, \alpha)$.
Furthermore, the range of a partial representation generates the same $C^*$-algebra as
the range of the corresponding covariant representation.
Since  $C^*_p(G; \mathcal R)$ and $C(\Omega_{\mathcal R}) \rtimes G$ are both 
generated by the ranges of universal representations, 
$C(\Omega_{\mathcal R}) \rtimes G$ is a realization of $C^*_p(G; \mathcal R)$.
\end{proof}

In the remaining sections we consider several situations that fall 
naturally into the framework of partial representations with relations.

\section{No relations: the partial group algebra $C^*_p(G)$.}
\label{no-relations}
Let $\mathcal R $ be the empty set of relations and
 consider {\em all} partial representations of a group $G$, 
subject to no restrictions. 
This is the situation considered in \cite{exe-pag} and
mentioned at the beginning of Section~\ref{spectrum-sect}.
The spectrum $\Omega_\emptyset$ is $X_G:= \{\omega \subset 2^G:
e \in \omega\}$ and the canonical partial action $\theta$ is given by 
$\theta_t (\omega) = t\omega$ for $\omega \ni t\inv$.
By \thmref{thm-spectrum}, the  crossed product $C(X_G) \rtimes_\alpha G$
is the  universal $C^*$-algebra $C^*_p(G) := C^*_p(G;\emptyset)$ 
for partial representations of $G$. 

\begin{proposition}
 The canonical partial action of a group $G$ on $X_G$ is topologically free 
if and only if $G$ is infinite.
\end{proposition}

\begin{proof}
When the group $G$ is finite, the spectrum $X_G$ has the 
discrete topology. Since the point $G \in X_G$ is
fixed by every group element, the partial action 
associated to partial representations of a finite group is never
topologically free.

Assume now that $G$ is infinite and let 
$$
U := \{ \chi \in X_G: a_i \in \chi \text{ and } b_j \notin \chi \text{ for }
1 \leq i \leq m \text{ and } 1 \leq j \leq n\}
$$
be a typical basic (nonempty) open set in $X_G$ where $a_i, b_j \in G$.
It suffices to show that for every element 
$t \in G\setminus\{e\}$ there is
some $\omega_0 \in U$ which is not fixed by $t$.
We may restrict our attention to the intersection of $U$ with the domain of
$\theta_t$, by assuming that one of the
$a_i$'s (and none of the $b_j$'s) is equal to $t\inv$.

Since $G$ is infinite there exists an element $c \in G$ different from
the $a_i$ and the $b_j$ and such that $t c $ is different from the $a_i$.
Then  $\omega_0 : = \{e, a_1, a_2, \cdots , a_m, c\}$ is in $U$
and $\theta_t(\omega_0) = \{t, ta_1, ta_2, \cdots , ta_m, tc\}$ is different from 
$\omega_0$ because $tc$ is not in $\omega_0$.
 Thus $\omega_0$ is not fixed by $t$, finishing the proof.
\end{proof}

\begin{corollary}
For infinite $G$,  a representation of $ C(X_G) \rtimes_r G$ is faithful
if and only if its restriction to $C(X_G)$ is faithful.
\end{corollary}
\begin{proof}
Direct application of  \thmref{thm-fax}.
\end{proof}
\begin{remark}
The singleton $\{G\} \subset  X_G$ is always closed and 
invariant under the partial action, so topological freeness fails 
at least for the restriction to $\{G\}$. Because of this 
the situation of \thmref{ideals}   
 never arises for the empty set of relations, and a characterization
of the ideals in $C_p^*(G)$ lies beyond the present techniques.
\end{remark}

\begin{theorem}
The canonical partial action of $G$ on $C(X_G)$ has the approximation property 
if and only if $G$ is amenable.
\end{theorem}
\begin{proof}
That the partial action of an amenable group $G$ 
on $X_G$ satisfies the approximation property is an easy consequence of 
\cite[Theorem 4.7]{exe-afb}.

To prove the converse suppose the action of $G$ on $X_G$ satisfies the 
approximation property. Then by  \proref{app-prop-quo}
the (trivial) action on the closed invariant singleton $\{ G\}$ 
satisfies the approximation property. Hence this trivial action
is amenable and the reduced and full crossed products coincide.
Since they correspond to the reduced and full group 
$C^*$-algebras of $G$, $G$ itself must be an amenable group.
\end{proof}

\begin{remark}
Since we do not know whether amenability itself is inherited by quotients, 
we do not know whether amenability of the partial action of $G$ on  $C(X_G)$
entails amenability of $G$. 
\end{remark}

\section{Nica covariance: the Toeplitz algebras of 
quasi-lattice groups.} 
\label{nica-rels}
Let $(G,P)$ be a quasi-lattice ordered group, as defined by 
Nica in \cite[\S 2]{nica}. The semigroup $P$ induces a partial order in $G$ 
via $x\leq y$ if and only if $x\inv y \in P$.  The quasi-lattice condition
says that if for $x, y \in G$ 
the set $\{z\in P: x \leq z, y \leq z\}$ is nonempty, then 
it has a unique minimal element, denoted $x \vee y$, 
and referred to as the least common upper bound in $P$ of $x$ and $y$,
(if there is no common upper bound, we write $x \vee y = \infty$). It is 
easy to see that $x$ has an upper bound in $P$ if and only if $x \in PP\inv$.

An {\em isometric representation} of $P$ on $H$ is a map $V: P \to B(H)$ 
such that $V_x^* V_x = 1$ and $V_x V_y = V_{xy}$.
 The isometric representation $V$ is {\em covariant} if it satisfies
$$
V_x V_x^* V_y V_y^* = V_{x\vee y} V_{x\vee y}^* \qquad x,y \in P,
$$
here we use the convention that $V_\infty = 0$, 
so that if $x$ and $y$ do not have a common upper bound in $P$ then
the corresponding isometries have orthogonal ranges.
  
The Toeplitz (or Wiener-Hopf) algebra ${\mathcal T}(G,P)$ is the $C^*$-algebra 
generated by the left regular representation 
$T$ of $P$ on $\ell^2(P)$ \cite{nica}, which is easily seen to be  covariant. 
The universal $C^*$-algebra $C^*(G,P)$ is the $C^*$-algebra generated
by a universal covariant semigroup of isometries.
When $(G,P)$ is amenable, the
canonical homomorphism $C^*(G,P) \mapsto {\mathcal T}(G,P)$ is 
an isomorphism \cite{nica,quasilat}.

Every $x \in PP\inv$ can be written in a ``most efficient way'' as
$x = \sigma(x) \tau(x)\inv$, where $\sigma(x) := x \vee e$ is the smallest
 upper bound of $x$ in $P$ and $\tau(x) := \sigma (x\inv) = x\inv \sigma(x)$.
 Using this factorization 
Raeburn and the third author have shown in \cite[Theorem 6.6]{qui-rae}
that ${\mathcal T}(G,P)$ is 
a crossed product  by a partial action on its diagonal subalgebra.
Their proof involves extending isometric covariant representations of $P$ to
 partial representations of $G$, and can be pushed further 
to describe the class of such extensions in terms of relations
satisfied by the range projections, which we do next. 

\begin{proposition} \label{pro-efficient}
Let $(G,P)$ be a quasi-lattice ordered group.
\begin{enumerate}
\item[(1.)] If $V$ is a covariant isometric representation of $P$ then 
\begin{equation}\label{efficient}
u_x = \left\{ 
\begin{array}{ll}
       V_{\sigma(x)} V_{\tau(x)}^*     &  \text{ if } x \in PP\inv \\
                             0         &  \text{ if } x \notin PP\inv.
\end{array}\right.
\end{equation}
is a partial representation of $G$ satisfying the relations
\begin{itemize}
\smallskip
\item[\textup($\mathcal N_1$\textup)]
  $\quad u_t^* u_t = 1$   for $t\in P$,   and    
\smallskip
\item[\textup($\mathcal N_2$\textup)] 
   $\quad u_x u_x^* u_y u_y^* = u_{x\vee y} u_{x\vee y}^*$   for $
x,y\in G$,
\smallskip
\end{itemize}
which we denote collectively by 
\textup($\mathcal N$\textup).
\item[(2.)] Conversely, every partial representation $u_t$ of $G$ 
satisfying the relations \textup($\mathcal N$\textup)
arises this way from a covariant isometric representation of $P$.
\end{enumerate}
\end{proposition}
\begin{proof}
(1.) That $u_x$ is a partial representation is proved in
 \cite[Theorem 6.6]{qui-rae} and that it satisfies ($\mathcal N_1$) is obvious.
We prove ($\mathcal N_2$) next.  Let $x,y \in G$ and assume both are in $P P\inv$, for
otherwise both sides are zero and there is nothing to prove. Notice first that
$u_x u_x^*  = V_{\sigma(x)} V_{\tau(x)}^* V_{\tau(x)} V_{\sigma(x)}^* = 
V_{\sigma(x)}  V_{\sigma(x)}^* $, so
$$
 u_x u_x^* u_y u_y^*  = V_{\sigma(x)}  V_{\sigma(x)}^* V_{\sigma(y)}  V_{\sigma(y)}^*
= V_{\sigma(x) \vee \sigma(y)}  V_{\sigma(x) \vee \sigma(y)}^* .
$$
Since $ x\vee y = x \vee e \vee y = \sigma(x) \vee \sigma(y)$ this proves ($\mathcal N_2$).

(2.) Assume now that $u_x$ is a partial representation of $G$ satisfying
 ($\mathcal N$). Then  $ u_t$ is an isometry for every $t\in P$,
and property (iii) in the 
Definition~\ref{def-p-rep} of partial representation gives 
$u_s u_t = u_s u_t u_t^* u_t = u_{st} u_t^* u_t = u_{st}$. 
Thus the restriction of $u$ 
to $P$ is an isometric representation, which is covariant by ($\mathcal N_2$).

It only remains to check that $u$ arises from its restriction $V$ to
$P$ as in (\ref{efficient}). Let $x \in G$. Then 
$u_x u_x^* = u_x u_x^* u_e u_e^* = u_{\sigma(x)} u_{\sigma(x)}^*$ by 
($\mathcal N_2$).  If $x \notin P P\inv$, then $\sigma(x) = \infty$ and   
 $u_x u_x^*$ vanishes. If $x \in P P\inv$, then 
$$
u_x = u_x u_x^* u_x = u_{\sigma(x)} u_{\sigma(x)}^* u_x. 
$$

The last two factors can be combined using
 property (iii) of Definition~\ref{def-p-rep}, 
and since $\sigma(x)\inv x = \tau(x) \inv$ we conclude that
$u_x = u_{\sigma(x)} u_{\tau(x)}^*$.
\end{proof}

\begin{definition}
A subset $\omega$ of $G$ is {\em hereditary}
if $x P\inv \subset \omega$  for every $x \in \omega$. It is 
{\em directed} if for every $x,y \in \omega$ there exists $z\in \omega \cap P$
with $x \leq z $ and  $y \leq z$.
\end{definition}

Notice that a hereditary subset $\omega$ is directed if and only if
the least upper bound of any two of its elements exists
and is in $\omega$; in particular, hereditary, directed 
subsets are contained in the set $P P\inv$.

\begin{lemma} The set of
hereditary, directed subsets of $G$ containing $e$
is invariant under the partial action $\theta$ on $X_G$. 
\end{lemma}
\begin{proof}
Suppose $\omega \in X_G$ is hereditary and directed and let $z\inv \in \omega$.
In order to see that $z\omega$ is hereditary, suppose $zx \in z\omega$ 
with $x \in \omega$ 
and let $t\in P$. Then $xt\inv \in \omega$ and   $zxt\inv \in z\omega$.

Next we show that $z\omega$ is directed. Assume $zx$ and $zy$ are
elements of $z\omega$.
Since $\omega$ is directed and contains $x$, $y$, and $z\inv$, 
it follows that $(x\vee y \vee z\inv) \in P\cap \omega$. 
It is easy to see using the definition that $z(x\vee y \vee z\inv) \in P \cap
z\omega$ is a common upper bound for $zx$ and $zy$. 
Thus $zx \vee zy \leq z(x\vee y \vee z\inv)$ and, since $z\omega $ is hereditary,
$zx \vee zy \in z\omega$.
\end{proof}

\begin{theorem}
The spectrum $\Omega_{\mathcal N}$ of the relations 
($\mathcal N$)
is the set
of hereditary, directed subsets of $G$ which contain the identity element. 

The  crossed product 
$C(\Omega_{\mathcal N}) \rtimes_\alpha G$ is canonically isomorphic to
the universal $C^*$-algebra $C^*(G,P)$ for covariant 
isometric representations of $P$.
\end{theorem}
\begin{proof}
Let $H$ be the set of hereditary, directed subsets of $G$ containing 
the identity element. Then clearly $H \subset X_G$. 

First we show that every $\omega \in \Omega_{\mathcal N}$ 
is  hereditary and directed;
 that $e\in \omega$ is obvious
because $\Omega_{\mathcal N} \subset X_G$.
If  $x \in \omega$ then $\omega$ is in the domain of the partial homeomorphism 
$\theta_{x\inv}$,
and since $\Omega_{\mathcal N}$ is invariant we have
$x\inv \omega = \theta_{x\inv}(\omega) \in \Omega_{\mathcal N}$.
By the relation 
($\mathcal N_1$), 
for $t \in P$ we obtain
 $[t]^*[t] (x\inv \omega) = 1$, which means that
$t\inv  \in x\inv \omega$. 
Since $xt\inv \in \omega$ for every $t\in P$ and every $x \in \omega$, 
$\omega$ is hereditary.

If $x$ and $y$ are elements of $\omega$,
then $1 = [x][x]^* [y][y]^* (\omega) = [x\vee y][x\vee y]^* (\omega)$
by ($\mathcal N_2$). Thus $x\vee y \in \omega$  and $\omega$ is directed.
 
Conversely, by the preceding lemma, 
if $\omega\in X_G$ is hereditary and directed, and if 
$z\inv \in \omega$, then $z \omega$ is also  hereditary and directed,
so it suffices to show that the relations ($\mathcal N$) hold
at every  hereditary, directed $\omega \in X_G$.

It is trivial to verify ($\mathcal N_1$), since $e\in \omega $ implies 
$et\inv \in \omega$ for every $t\in P$ by hereditarity of $\omega$.
For ($\mathcal N_2$) we need to show that
$ [x] [x]^* [y] [y]^* (\omega )  = [x\vee y] [x\vee y]^* (\omega)$, or,
equivalently, that $x$ and $y$ are in $\omega$
 if and only if $x\vee y \in \omega$. The ``only if''
holds because $\omega$ is directed,
and the ``if'' holds because it is hereditary,
since $(x\vee y)\inv x \in P$ and
$x = (x\vee y) \, (x\vee y)\inv x$. 

The  crossed product is isomorphic to
$C^*(G,P)$ because of \proref{pro-efficient}.
\end{proof}

\begin{remark} Hereditary directed subsets of the semigroup $P$ were introduced
by Nica in \cite[\S 6.2]{nica}, where he showed that the spectrum of
the diagonal subalgebra in the Toeplitz 
algebra
is (homeomorphic to) the
space of hereditary, directed, and nonempty subsets of $P$.
The homeomorphism of our spectrum $\Omega_{\mathcal N}$ 
to 
the space considered by Nica is obtained simply by  
sending an element $\omega$ of $\Omega_{\mathcal N}$ to its intersection with $P$.
\end{remark}

\begin{proposition}
The canonical partial action $\theta$ on $\Omega_{\mathcal N}$
is topologically free.
\end{proposition}
\begin{proof}
For each  $t \in P$ the set
 $tP\inv = \{x\in G: x \leq t\}$  is hereditary and directed;
moreover, $t\neq t'$ implies $tP\inv \neq t'P\inv$. This gives a copy of $P$
inside $\Omega_{\mathcal N}$ which is in fact 
dense
\cite[\S 6.2]{nica}.

Suppose $x \in G$; it is easy to see that the point $tP\inv $ is in the 
domain of the partial 
homeomorphism
$\theta_x$ if and only if 
$xt\in P$. In this case $\theta_x(tP\inv) = (xt) P\inv \neq tP\inv$.
Since no point in this dense subset is fixed by $\theta_x$ for 
$x\neq e$,
the proof is finished.
\end{proof}

As an application  
we obtain a characterization of faithful representations
of the reduced crossed product which 
is slightly more general than \cite[Theorem 3.7]{quasilat}, in 
that it does away with the amenability hypothesis by
focusing on the reduced crossed product.
{}From this point of view, it becomes apparent 
that the faithfulness theorem for representations 
is really a theorem about {\em reduced} crossed products
and that it is a manifestation of topological freeness.
\begin{theorem} 
Suppose  $(G,P)$ is a quasi-lattice ordered group. 
A representation of the reduced crossed product
 $C(\Omega_{\mathcal N}) {\rtimes}_{\alpha,r} G$ 
is faithful 
if and only if it is faithful on the diagonal $C(\Omega_{\mathcal N})$.
\end{theorem}
\begin{proof}
Since $\alpha$ is topologically free, the result follows from \thmref{thm-fax}.
\end{proof}
Of course we may use \cite[Proposition 2.3(3)]{quasilat} to 
decide when the restriction to $C(\Omega_{\mathcal N})$ 
of a representation $\pi \times v$ of
$C(\Omega_{\mathcal N}) {\rtimes}_{\alpha,r} G$ is faithful in terms of the 
generating partial isometries: the condition is that
$$
\pi ( {\textstyle\prod_{t\in F} } ( 1 - u_t u_t^*) ) \neq 0
$$
for every finite subset $F$ of $P \setminus \{e\}$.

Since the diagonal algebra in Nica's Wiener-Hopf $C^*$-algebra $T (G,P)$
is a faithful copy of $ C(\Omega_{\mathcal N})$ 
we can use the faithfulness theorem to express $T (G,P)$
as the reduced crossed product by a partial action. 
\begin{corollary} 
If $(G,P)$ is a quasi-lattice ordered group, then
$$ C(\Omega_{\mathcal N}) {\rtimes}_{\alpha,r} G \cong \mathcal T (G,P).$$
\end{corollary}

This isomorphism is essentially \cite[Theorem~6.6]{qui-rae};
although the partial action there is not given explicitly,
it is not hard to see that it is the one above.

\section{Cuntz-Krieger relations: the universal ${\mathcal O}_A$.}
\label{ck-section} Let $A=[a_{ij}]$ be  a $\{0,1\}$-valued 
$n$ by $n$ matrix with no zero rows. A {\em Cuntz-Krieger $A$-family} is a
collection  of partial isometries $\{s_i\}_{i=1}^n$  such that
$$
\sum_j s_j s_j^* = 1, \quad \text{ and } \quad
\sum_j a_{ij} s_j s_j^* = s_i^* s_i \quad\text{ for }i = 1, 2, \ldots ,n.
$$

We define the algebra $\mathcal O_A$ to be the
universal $C^*$-algebra for Cuntz-Krieger $A$-families. 
That is, $\mathcal O_A$ is
the $C^*$-algebra  generated by $n$ partial isometries 
$\{s_i\}_{i=1}^n$ satisfying the two conditions above,
and such that if $\{s'_i\}_{i=1}^n$ is any other
collection satisfying the same conditions, the map $s_i \mapsto s'_i$ extends to 
a $C^*$-algebra homomorphism from $C^*( \{s_i\}_{i=1}^n)$  onto 
$C^*(\{s'_i\}_{i=1}^n)$. 

This is not the original
definition of $\mathcal O_A$, given in \cite{cun-kri} only for $A$ satisfying a
certain condition (I) which implies the uniqueness of the
$C^*$-algebra $C^*(\{s_i\})$ provided that $s_i \neq 0$ for every $i$. 
We have chosen to define $\mathcal O_A$ as a universal object,
as in \cite{anh-rae}, because it allows us to treat more general 
matrices,
and at the same time clarifies the presentation. With this 
definition, Cuntz and Krieger's celebrated uniqueness result
says that if $A$ satisfies (I) then  $\mathcal O_A$ is isomorphic to the
$C^*$-algebra generated by {\em any} $A$-family of nonzero partial isometries.

We will prove that $\mathcal O_A$ is a crossed product by a partial action of
the free group $\mathbb F_n$ on the space of infinite admissible paths, by
showing first that it is the universal $C^*$-algebra of partial
representations of $\mathbb F_n$ subject to certain relations, and then using
\thmref{thm-spectrum} to compute the spectrum and the partial action.  This
will allow us to obtain Cuntz and Krieger's uniqueness theorem from our
\thmref{thm-fax}.  First we need to set some notation and recall some
terminology.

\begin{definition}
Let $|x|$ denote the length of a reduced word $x \in {\mathbb F}_n$.
A partial representation $u$
of ${\mathbb F}_n$ is {\em semisaturated} if it satisfies
$$u_{tr} = u_t u_r \quad \text{ whenever } \quad |tr| = |t| + |r|,$$
and, similarly,  
if the partial action $\alpha$ of $\mathbb F_n$ satisfies
$$\alpha_{tr} = \alpha_t \circ \alpha_r 
\quad \text{ whenever } \quad |tr| = |t| + |r|,$$
we say that it is a {\em semisaturated partial action}.
The condition $|tr| = |t| + |r|$ means that there is no further reduction 
in the concatenation of the reduced words $t$ and $r$.
\end{definition}

Unlike general partial representations, semisaturated ones are determined by their
values on the generators $\{g_1,\dots,g_n\}$ of the free group: if 
$t = g_i^{\pm 1} g_j^{\pm 1} \cdots g_k^{\pm 1}$ is a reduced word in
$\mathbb F_n$ then $u_t = u_{g_i^{\pm 1}} u_{g_j^{\pm 1}} \cdots u_{g_k^{\pm 1}}$,
and similarly for a semisaturated partial action.
For this reason semisaturated partial actions are 
also called {\em multiplicative} in \cite{qui-rae}.

Semisaturated partial representations
and partial actions can be characterized in terms of the 
range projections $e(t) = u_{t} u_{t}^*$, as in 
 \cite[Proposition 5.4]{exe-afb} and \cite[Lemma 5.2]{qui-rae}: 
$u$ is semisaturated if and only if
$e(tr) e(t) = e(tr)$ whenever $ |tr| = |t| + |r|$,
and similarly for partial actions.

\begin{proposition}{\em (cf. \cite[Theorem 5.2]{exe-afb})}\label{ck-part-rep}
 Let $\{ g_1, g_2, \ldots g_n\}$ be the generators of $\mathbb F_n$.
For every  Cuntz-Krieger 
$A$-family $\{s_i\}_{i=1}^n$ there is a  unique semisaturated
partial representation $u$ of ${\mathbb F}_n$ such that $u_{g_i} = s_i$.
The range projections $e(g) = u_g u_g^*$  satisfy the relations 
\begin{equation}
e(tr) e(t) = e(tr) \qquad \text{ whenever } |tr| = |t| + |r|,
 \tag{$\mathcal {CK}_{ss}$}
\end{equation}
\begin{equation}
 \sum_{j=1}^n e(g_j) = 1, \qquad \qquad \qquad \qquad and
\tag{$\mathcal {CK}_{1}$}
\end{equation}
\begin{equation}
 \sum_{j=1}^n a_{ij}  e(g_j)  = e(g_i\inv) 
\qquad \quad \text{ for } i = 1, \ldots , n,
\tag{$\mathcal {CK}_{A}$}
\end{equation}
which we denote
collectively by \textup($\mathcal {CK}$\textup).

Conversely, if $u$ is a partial representation of ${\mathbb F}_n$ satisfying these
relations then the partial isometries $s_i := u_{g_i}$ 
form a Cuntz-Krieger $A$-family. 

Thus the Cuntz-Krieger algebra $\mathcal O_A$ 
\textup(i.e., the universal $C^*$-algebra of Cuntz-Krieger
$A$-families\textup)
is universal for partial representations of $\mathbb F_n$
subject to the  relations \textup($\mathcal {CK}$\textup). 
\end{proposition}

 From this, \thmref{thm-spectrum} implies that 
there is a partial action $\theta$ on the spectrum 
$\Omega_{\mathcal {CK}}$ such that $s_i \mapsto u_{g_i}$
extends to an isomorphism of $\mathcal O_A$ to the crossed product 
$C(\Omega_{\mathcal {CK}}) \rtimes {\mathbb F}_n$. 
Clearly this partial action is semisaturated
because of ($\mathcal {CK}_{ss}$) and \cite[Lemma 5.2]{qui-rae}.

Although a direct application of Proposition \ref{spec-prop} naturally  yields
the spectrum  $\Omega_{\mathcal {CK}}$ of the relations ($\mathcal {CK}$) as a
subspace of $X_{\mathbb F_n}$, we will eventually identify this spectrum with the
more familiar infinite path space, so we will state the result in terms of path
space, which we define next.

\begin{definition}
An {\em infinite \textup(admissible\textup) path} 
is an infinite sequence $\mu = \mu_1 \mu_2 \ldots $ of generators
of $\mathbb F_n$ such that $a(\mu_j,\mu_{j+1}) = 1$ for every $j\in\mathbb N$.
{\em Infinite path space} is the space $P_A^\infty$ of infinite admissible paths
with the relative topology inherited as a closed and hence compact
subspace of the infinite product space $\prod_{i=1}^\infty \{g_1, \ldots, g_n\}$.
\end{definition}

\begin{theorem}\label{CKalg=pcp}
There is a unique semisaturated partial action $\theta$ of the free group 
$\mathbb F_n$ on infinite path space $P_A^\infty$ such that 
\begin{equation}\label{CK-phom}
\theta_{g_i}(\mu) = g_i\mu  \qquad \text{ for } \mu \in U_{g_i\inv}:=
\{\mu: a(g_i,\mu_1) = 1\} \subset P_A^\infty,
\end{equation}
where $g_i \mu$ means concatenation of $g_i$ at the beginning of $\mu$.
The map $s_i \mapsto u_{g_i}$ extends to an isomorphism of $\mathcal O_A$ 
to the crossed product $C(P_A^\infty) \rtimes_\theta \mathbb F_n$.
\end{theorem}

The theorem will be proved in the following lemma and propositions
by establishing a homeomorphism of $\Omega_{\mathcal {CK}}$  to $P_A^\infty$,
and then showing that the partial action carried over to $P_A^\infty$,
which is necessarily semisaturated, 
satisfies (\ref{CK-phom}).

\begin{definition}
A subset $\omega$ of ${\mathbb F_n}$ is {\em connected} if, when viewed as a
subset of the Cayley graph of ${\mathbb F_n}$, it contains the shortest path 
between any two of its elements.
\end{definition}

When $e \in \omega$ it is easy to see that $\omega$
is connected if and only if it contains the
initial subwords of the reduced words in $\omega$.
Next we show that the semisaturation relations single out the
{\em connected} subsets   in $X_{\mathbb F_n}$.
\begin{lemma}\label{ss-spectrum}
 The relations 
$$
e(t) e(tr) = e(tr) \quad \text{ for }
t,r \in \mathbb F_n  \text{ such that } |tr| = |t| + |r|
$$ 
are satisfied at the point $\omega \in X_{\mathbb F_n}$
if and only if $\omega$ is connected as a subset of ${\mathbb F}_n$.
\end{lemma}

\begin{proof}
Suppose  $|tr| = |t| + |r|$. Then 
\begin{equation*}
(1_{tr} 1_t) |_\omega = 1_{tr}|_\omega 
\quad  \iff \quad 
(tr \in \omega  \Longrightarrow t\in \omega).
\end{equation*}
Since each $\omega \in X_{\mathbb F_n}$ contains $e$, it is closed under 
taking initial subwords if and only if it is connected. 
\end{proof}

\begin{proposition}\label{ck-spectrum}
A subset $\omega \subset {\mathbb F}_n$ is in the spectrum 
$\Omega_{\mathcal {CK}}$ of the Cuntz-Krieger relations if and 
only if the following conditions hold
\begin{itemize} 
\item[\spec{1}]  $ e \in \omega$,

\item[\spec{2}]  $\omega$ is connected,

\item[\spec{3}]  for every $t\in \omega$  there is a unique generator 
      $g_j =  g_{j(\omega,t)}$ such that $tg_j \in \omega$, and

\item[\spec{4}]  for every $ t\in \omega $, we have that 
        $ tg_i\inv \in \omega$ if and only if $ a_{i,{j(\omega,t)}}  = 1 $,  with
        ${j(\omega,t)}$ given in {\em \spec{3}}. 

\end{itemize}

\end{proposition}
\begin{proof}
Let $Z$ be the set of $\omega$'s satisfying conditions
\spec{1}{}--\spec{4}; we need to show $Z = \Omega_{\mathcal {CK}}$.
The first step is to show that $Z$ is $\theta$-invariant.
To do this, suppose $\omega$ satisfies \spec{1}--\spec{4},
 and let  $x\inv \in \omega$. Then
$\theta_x(\omega) = x\omega$ also satisfies
\spec{1} because $x\inv \in \omega$; \spec{2} also holds for $x\omega$
because connected sets are  translation-invariant.
It is routine to check that $j(\omega,t)$ depends only on
$t\inv \omega$, hence  $j(x\omega,xt) = j(\omega,t)$, from which 
\spec{3} and \spec{4} with $x\omega$ in place of $\omega$ follow easily .

Assume now $\omega \in Z$.  Condition \spec{2} implies  that the
semisaturation relation holds at $\omega$, from \lemref{ss-spectrum}. 
{}From \spec{3} with $t=e $, there exists a unique $j(\omega,e)$ such that
$g_{j(\omega,e)} \in \omega$, that is, such that $1_{g_{j(\omega,e)}} |_\omega = 1$,
{}from which ($\mathcal {CK}_1$) holds at $\omega$.
Setting $t = e$ in \spec{4}, $g_i\inv \in \omega $ if and only if 
$a_{i,{j(\omega,e)}} = 1$, from which ($\mathcal {CK}_A$)  holds at $\omega$.
This implies that $ Z \subset \Omega_{\mathcal {CK}}$,
because it is an invariant closed subset on which the relations are satisfied. 

If $\omega \in \Omega_{\mathcal {CK}}$ then 
$\sum_j 1_{g_j} (t\inv \omega)  = 1$ for $t\in \omega$
by ($\mathcal {CK}_1$) and (\ref{spec-def}), so there is exactly one
generator $g_j$ such that $1_{g_j} (t\inv \omega)  = 1$, i.e., 
such that $tg_j \in \omega$, proving \spec{3}.
{}From ($\mathcal {CK}_A$) and (\ref{spec-def}),
$$
 \sum_{j=1}^n a_{ij} 1_{g_j} (t\inv \omega) = 1_{g_i\inv}(t\inv \omega),
$$ and since  $1_{g_j}( t\inv \omega) \neq 0$ only for 
$j = j(\omega,t)$, \spec{4} follows.
\end{proof}

\begin{proposition}\label{no-Lcorners}
Suppose $y \in \mathbb F_n$. If the partial homeomorphism $\theta_{y\inv}$ of 
$\Omega_{\mathcal {CK}}$ is nontrivial then $y = r s \inv$ for 
two admissible words $r,s \in \mathbb F_n^+$ with the same final letter. 
\end{proposition}
\begin{proof}
If the partial automorphism $\theta_{y\inv}$ is nontrivial then  
$y  \in\omega$ for some $ \omega \in \Omega_{\mathcal {CK}}$. 
Suppose somewhere in the reduced form of $y$ there is an inverse generator
 followed by a generator so that $y$ has an initial reduced subword 
of the form $x g\inv g'$. By connectedness of $\omega$ the subwords $x$ and
$xg\inv$ are also in $\omega$. But this violates condition \spec{3} 
at $t = xg\inv$: the ``forward continuations'' $x = (x g\inv) g$ and 
$x_1 = (xg\inv)g'$ are different because $g \neq g'$ and they are both in $\omega$.
Hence  the reduced form of $y$ is $r s\inv$ with $r$ and $s$ admissible.
By connectedness $r\in \omega$, and by \spec{3} there is a  unique generator $g$
such that $rg\in \omega$. Let $s_\ell $ be the last letter of 
$s$; again by connectedness, $rs_\ell\inv \in \omega$, and it follows from \spec{4}
that $a(s_\ell,g) = 1$. This makes $sg$ an admissible word, and since
$y = (rg) (sg)\inv$ the proof is finished.
\end{proof}

We aim to show that $\Omega_{\mathcal {CK}}$ is covariantly homeomorphic 
to infinite path space with the partial action given in (\ref{CK-phom}).
The appropriate homeomorphism is defined as follows.

Since  every $\omega \in \Omega_{\mathcal C \mathcal K}$ contains the identity
element and satisfies \spec{3}, it contains a unique generator $g_{i_1}$, 
and, again, a unique product $g_{i_1} g_{i_2}$, etc.
By induction, for each $\omega \in \Omega_{\mathcal C \mathcal K}$ there exists a
unique  infinite admissible path $\omega^+$ such that every finite initial segment
of $\omega^+$,  when viewed as an element of $\mathbb F_n$, belongs to $\omega$.

We may identify the infinite path $\mu = \mu_1  \mu_2   \mu_3  \cdots $
with the collection of its finite initial segments $\{ \mu_1,\, \mu_1 \mu_2,\, 
\mu_1 \mu_2 \mu_3,\, \ldots \}$, and with this picture in mind the correspondence
$\omega \mapsto \omega^+$ mentioned above comes simply from intersecting
$\omega\in \Omega_{\mathcal {CK}}$ with  $\mathbb F_n^+$, 
i.e., $\omega^+ = \omega \cap \mathbb F_n^+$.

\begin{proposition}
The map $\omega \mapsto \omega \cap {\mathbb F_n^+} $ is a homeomorphism
of $\Omega_{\mathcal {CK}}$ onto $P_A^\infty$; the corresponding partial action 
 $\tilde{\theta}$ on $P_A^\infty$ is the unique
semisaturated one for which a generator $g$ acts by
$\tilde{\theta}_{g }(\mu) = g \mu$ on the set 
$U_{g \inv} := \{\mu: a(g,\mu_1) = 1\} \subset P_A^\infty$.
\end{proposition}
\begin{proof}
We show first that the map $\omega \mapsto \omega \cap \mathbb F_n^+$ is injective.
Assume that $\omega $ and $\omega'$ are in $\Omega_{\mathcal {CK}}$ and 
$\omega \cap \mathbb F_n^+ = \omega' \cap \mathbb F_n^+$. We will show that 
$\omega \subset \omega '$ and hence, by symmetry, that they coincide.
Let $t \in \omega$. From  \proref{no-Lcorners}
we know that $t = \mu \nu\inv$ for admissible finite paths $\mu, \nu$ 
in $\mathbb F_n^+$. If $|\nu| = 0$, then $t = \mu \in \mathbb F_n^+$, 
so $t \in \omega '$. If $\nu$ is not the empty word, let $\nu_\ell$
be its last letter.
Since $\omega$ is connected, $\mu$ is in $\omega$, and since $\mu$ is positive,
$\mu \in \omega \cap \mathbb F_n^+ = \omega' \cap \mathbb F_n^+$, so
$\mu \in \omega'$. It follows from \spec{3}
that the next admissible generator after $\mu$ 
is the same in both $\omega$ and $\omega'$, i.e., 
$j(\mu,\omega) = j(\mu,\omega')$. Since $\mu \nu_\ell\inv$ is in $\omega$ by
connectedness (because $t = \mu \nu\inv \in \omega$), 
we have  $a(\nu_\ell, j(\omega,\mu,)) = 1$, so
condition \spec{4} puts $\mu \nu_\ell\inv $ in $\omega'$.
Continuing the argument by induction gives that
$t = \mu \nu_\ell\inv \ldots \nu_1\inv$ is in $\omega '$, as claimed.

In order to show that the map $\omega\mapsto\omega\cap\mathbb F_n^+$
is onto, suppose $\mu$ is an infinite admissible path, and let 
$$
\omega_\mu := \{ t \in \mathbb F_n : t = (\mu_1 \mu_2 \cdots \mu_k) (\nu)\inv  
\text{ with } k \geq 0 \text{ and } \nu \mu_{k+1} \text{ admissible }\}.
$$ 
It is easy to show that $\omega_\mu$ satisfies conditions  
\spec{1}--\spec{4}, so it is in $\Omega_{\mathcal{CK}}$, and clearly
$\omega_\mu \cap \mathbb F_n^+ = \mu$.

Let $\mu \in P_A^\infty$ and $k \geq 1$. The basic neighborhood 
$$
V_{\mu,k}: =  \{\nu\in P_A^\infty: \nu_i = \mu_i \text{ for }i\leq k\}
$$
of $\mu$ is the image of the set 
$\{\omega \in \Omega_{\mathcal {CK}}: \mu_1 \mu_2 \cdots \mu_k \in
\omega\}$, which is clearly open in $\Omega_{\mathcal {CK}}$.
Hence the map $\omega \mapsto \omega^+$ is continuous, and since it is 
a bijection from a compact to a Hausdorff space, it is a homeomorphism.

 The corresponding partial action $\tilde \theta$ on $P_A^\infty$ 
is semisaturated, and hence is characterized by its behaviour on generators.
If $g$ is a generator of $\mathbb F_n$
then $\tilde \theta_g (\mu) = \theta_g(\omega_\mu)^+ 
= (g\omega_\mu)^+ =  g\mu$ for every $\mu$ such that $a(g, \mu_1) =1$. 
\end{proof}

\begin{proof} 
(of \thmref{CKalg=pcp})
By \proref{ck-part-rep}, $\mathcal O_A$ is 
isomorphic to the  crossed product 
$C(\Omega_{\mathcal {CK}}) \rtimes  \mathbb F_n$,
and by the preceding proposition we may replace 
$\Omega_{\mathcal {CK}}$ by $P_A^\infty$
with the partial action defined in (\ref{CK-phom}).
\end{proof}

Next we discuss some basic properties of paths with an eye 
to characterizing the matrices for which the
partial action is topologically free.

\begin{definition}
An infinite path $\mu \in P_A^\infty$ is {\em periodic}
if it is of the form $\mu = x \gamma \gamma \gamma \cdots$
for $x$ and $\gamma$ finite admissible words.
A path is  {\em aperiodic} if it is not periodic.
\end{definition}
Notice that the word $\gamma$ above is an admissible {\em circuit }
(this means $\gamma \neq e$ and $\gamma \gamma$ is admissible).
The words $x$ and $\gamma$ in the definition of periodic paths
 are not unique, as evidenced by 
placing some parentheses: e.g., $x \gamma \gamma \gamma \cdots  
= (x\gamma) \gamma^2 \gamma^2 \cdots $.
However, one may reduce a periodic path in a canonical way
by requiring that $\gamma$ be the shortest possible and that $x$
contain  no part of $\gamma$ in the sense that the last letters of 
$x$ and $\gamma$ are different --- if they are not then shorten $x$ and shift back
the origin of $\gamma$. (The same effect is obtained by reducing the product
$x \gamma x\inv$.)
 
\begin{definition}
An admissible circuit $\gamma = \gamma_1 \gamma_2 \cdots  \gamma_k$
is {\em terminal} if the row corresponding to each $\gamma_i$ in the matrix 
$A$ has only one nonzero entry. 
\end{definition} 
Note that, in the graph associated to $A$,
 a terminal circuit is a circuit with no exit.

\begin{lemma}\label{noexit=isolated}
The infinite admissible path $\mu$ is an isolated point in $P_A^\infty$ if and only if
$\mu = x \gamma \gamma \gamma \cdots $ for some terminal circuit $\gamma $ 
and some admissible word $x$. 
\end{lemma}
 
\begin{proof}
If $\mu = x \gamma \gamma \gamma \cdots \in P_A^\infty$ with
 $\gamma$ an admissible terminal circuit then $\mu$ is an isolated point in 
 $P_A^\infty$, because it is the only admissible path in the open set 
 $U := \{ \nu \in P_A^\infty: x\gamma \text{ is an initial segment of }\nu\}$. 
  
Suppose  $\mu \in P_A^\infty$ is aperiodic. Since $\mu$ is an infinite path, 
it will describe a circuit at least every $n$ steps because there are only
$n$ generators. Because of this, the initial segments of $\mu$  that
can be continued by repeating  a circuit 
are of unbounded length. Hence every neighborhood of $\mu$ contains 
a periodic (admissible, infinite) path which is necessarily
different from $\mu$, which was assumed aperiodic.
Therefore $\mu$ is not isolated.  
\end{proof}

Periodic paths are also the only possible nontrivial 
fixed points of the partial action:

\begin{lemma}\label{fix-point}
The path $\mu \in P_A^\infty$ is fixed by a nontrivial element 
$t\in \mathbb F_n$ under the partial action $\theta$ 
if and only if there exist admissible words $x$ and $\gamma$ such that 
$t = x\gamma^k x\inv$ with $k = \pm 1$ 
and $\mu = x \gamma \gamma \gamma \cdots$.
\end{lemma}

\begin{proof}
Assume $\mu$ is a periodic path in reduced form $x \gamma \gamma \gamma \cdots  $
and $t = x \gamma^k x\inv$ 
(we may suppose 
$k = 1$, otherwise change $t$ to $t\inv$).
Then $\theta_{x\inv}(\mu) = \gamma \gamma \gamma \cdots  $
and $\theta_{x\gamma^k} (\gamma \gamma \gamma \cdots) = x \gamma^k \gamma \gamma \gamma
\cdots = \mu$. Since $\theta_{x \gamma^k x\inv} $ extends 
$\theta_{x\gamma^k} \circ \theta_{x\inv}$ we have $\theta_t (\mu ) = \mu$.

To prove the converse, suppose $\theta_t$ fixes the infinite path $\mu$. 
By \proref{no-Lcorners}, $t = r s\inv$ (in reduced form) 
for two admissible words $r,s$.
Since $\mu$ is in the domain of $\theta_t$, 
 $s$ must be an initial segment of $\mu$, i.e., there exists $\mu'$ such that
  $\mu = s \mu'$ and hence
 $rs\inv \mu = r \mu'$. Thus  $r\mu' = s \mu'$,
so $r$ and $s$ have the same first letter, and second, etc., and
one of the two must run out of letters before the other, for otherwise they
would  coincide.
Say that $r$ is longer, then $r = s \gamma$ and we have
$s \gamma \mu' = s \mu'$, hence $\gamma \mu' = \mu'$ with 
$\gamma \in \mathbb F_n^+$, which
by induction  implies that $\mu' = \gamma \gamma \gamma \cdots$.
Hence $\mu = s \gamma \gamma \gamma \cdots$ and $t = s\gamma s\inv$.
\end{proof}

Cuntz and Krieger's condition (I) is known to be equivalent to
density of the aperiodic  paths, to the absence of isolated points, 
and to the absence of terminal circuits for the matrix $A$.
In the following proposition  
topological freeness of the partial action is added to the list.
\begin{proposition}\label{topfree-iff}
The following are equivalent:
\begin{enumerate}
\item The partial action $(C(P_A^\infty), \mathbb F_n, \theta)$ is 
       topologically free.
\item The graph with incidence matrix $A$ has no terminal circuits.
\item There are no isolated points in $P_A^\infty$.
\item The aperiodic paths are dense in $P_A^\infty$.
\item The matrix $A$ satisfies Cuntz and Krieger's condition {\em(I)}.
\end{enumerate}
\end{proposition}
\begin{proof}
By \lemref{fix-point}, a fixed point of the partial action is
determined by the group element that fixes it. Specifically,
the group element $t $ has a nonempty fixed point set if and only if
$t = x \gamma x\inv $ with $x\gamma \gamma$ admissible.
When this happens, the fixed point set for $t$
is just the singleton $\{ x  \gamma \gamma \gamma \cdots\} \subset P_A^\infty$.
By \lemref{noexit=isolated} this singleton is open if and only if $\gamma$ 
is terminal, so (i) $\iff $ (ii) $ \iff $ (iii). 

The proof of  (ii) $ \iff $ (v) can be found in
\cite[Lemma~3.3]{kpr} and  (iv) $ \iff $ (v) is from \cite{cun-kri}.
\end{proof}

As the main result of this section we 
obtain the Cuntz-Krieger uniqueness theorem 
via  an application of our characterization of 
faithful representations of crossed products
of topologically free partial actions.  

\begin{theorem}{\em (\cite[Theorem 2.13]{cun-kri})}
Suppose $A$ has no terminal circuits.
If $\{s_i\}_{i=1}^n$ and $\{s'_i\}_{i=1}^n$ are two Cuntz-Krieger $A$-families 
of nonzero partial isometries, then 
the map $s_i \mapsto s'_i$ gives an isomorphism of 
$C^*(\{s_i\}_{i=1}^n)$ to $C^*(\{s_i'\}_{i=1}^n)$.
\end{theorem}
\begin{proof}
The canonical partial action $\theta$ on $P_A^\infty$
is topologically free by \proref{topfree-iff}, and the reduced and full crossed
products coincide by \cite[Theorem 6.6]{exe-afb}. 
 
It suffices to show that nonzero partial isometries
give rise to faithful representations of $C(P_A^\infty)$, because then 
\thmref{thm-fax} implies that the representations 
of $C(P_A^\infty) \rtimes \mathbb F_n$
arising from $\{s_i\}$ and $\{s_i'\}$ are both faithful.

Everything hinges upon showing that if $(\pi,u)$ is a covariant representation
in which $\pi$ is not faithful, then $u_g =0 $ for some generator $g$ of $\mathbb
F_n$. Suppose $\pi$ is not faithful. Then there exists 
$f\in C(\Omega_{\mathcal {CK}})$ such that $f\neq 0$ and $\pi(f) = 0$.  
Without loss of generality we may assume $f$ to be positive.
Let $\mu \in P_A^\infty$ be a point with $f(\mu) > \|f\|/2 $. 
Then there exists $k$ large enough such that
 $f > \|f\|/2$ on the neighborhood 
$$
V_{\mu,k}: =  \{\nu\in P_A^\infty: \nu_i = \mu_i \text{ for }i\leq k\}
$$
of $\mu$.  Let $s = \mu_1 \mu_2 \cdots \mu_k$. Then 
$0 \leq u_s u_s^* = \pi(1_s) \leq (2/\|f\|) \pi(f) = 0$, so
$ u_s^* u_s = 0$. Since $  u_s^* u_s = u_{\mu_k} ^* u_{\mu_k}$
by \cite[Lemma 2.1 (a)]{cun-kri}, we conclude
 that the partial isometry  $u{\mu_k}$ vanishes.
\end{proof}

\begin{remark} It is also not hard to see that if $A$ is irreducible
and not a permutation matrix, 
then the partial action on $P_A^\infty$ is  minimal and topologically free, 
so the simplicity result for $\mathcal O_A$ 
\cite[Theorem 2.14]{cun-kri} follows from our Corollary~\ref{simple}.
\end{remark}

\end{document}